\documentclass[a4paper,12pt]{article}
\usepackage{amsthm,amsfonts,amssymb,euscript}
\usepackage{latexsym, multicol, fancybox}
\usepackage{graphicx}
\usepackage{color}
\usepackage{amsmath, amsthm, amssymb, bm}
\usepackage{epstopdf}
\usepackage{caption}
\usepackage{psfrag}

\usepackage{mathrsfs}
 \usepackage{xcolor}
 \usepackage[citebordercolor={green}]{hyperref}
 \usepackage{tikz}
 \usepackage{bbm}
 \usepackage{dsfont}
 \usepackage{bbold}

\numberwithin{equation}{section}

\newtheorem{theorem}{Theorem}[section]

\newtheorem{definition}[theorem]{Definition}
\newtheorem{remark}[theorem]{Remark}

\newcommand{\bea}{\begin{eqnarray}}
\newcommand{\eea}{\end{eqnarray}}
\def\beaa{\begin{eqnarray*}}
\def\eeaa{\end{eqnarray*}}
\def\ba{\begin{array}}
\def\ea{\end{array}}
\def\be#1{\begin{equation} \label{#1}}
\def \eeq{\end{equation}}

\def\a{{\alpha}}

\def\b{{\beta}}
\def\be{{\beta}}
\def\ga{\gamma}
\def\Ga{\Gamma}
\def\de{\delta}
\def\De{\Delta}
\def\ep{\epsilon}

\def\la{\lambda}
\def\La{\Lambda}

\def\Si{\Sigma}
\def\om{\omega}

\def\vphi{\varphi}

\def\th{\theta}

\def\ze{\zeta}

\def\nab{\nabla}

\def\pr{{\partial}}

\def\c{\cdot}

\def\AA{{math\cal A}}
\def\BB{{\mathcal B}}
\def\CC{{\mathcal C}}
\def\MM{{\mathcal M}}
\def\NN{{\mathcal N}}

\def\LL{{\mathcal L}}
\def\II{{\mathcal I}}

\def\HH{{\mathcal H}}

\def\TT{{\mathcal T}}

\def\JJ{{\mathcal J}}
\def\KK{{\mathcal K}}
\def\Lie{{\mathcal L}}

\def\AA{{\mathcal A}}

\def\HH{{\mathcal H}}

\def\Lie{{\mathcal L}}

\def\C{{\bf C}}
\def\D{{\bf D}}

\def\G{{\bf G}}

\def\R{{\bf R}}

\def\T{{\bf T}}
\def\Z{{\bf Z}}

\def\g{{\bf g}}

\def\RRR{{\Bbb R}}

\def\f12{{\frac 1 2}}

\def\dual{{\,\,^*}}

\def\Hb{\,\underline{H}}

\def\Xh{\,^{(h)}X}

\def\trch{{\mbox tr}\, \chi}
\def\chih{{\widehat \chi}}
\def\chib{{\underline \chi}}
\def\chibh{{\underline{\chih}}}

\def\etab{{\underline \eta}}
\def\omb{{\underline{\om}}}
\def\bb{{\underline{\b}}}
\def\aa{\protect\underline{\a}}
\def\xib{{\underline \xi}}

\def\Ab{\protect\underline{A}}
\def\Bb{\protect\underline{B}}
\def\Xb{\protect\underline{X}}

\newcommand{\CCb}{\underline{\mathcal{C}}}
\newcommand{\BBb}{\underline{\mathcal{\BB}}}

\def\Xh{\widehat{X}}
\def\Xbh{\widehat{\Xb}}

\def\ub{\underline{u}}

\def\tr{\mbox{tr}}
\def\atr{\,^{(a)}\mbox{tr}}

\def\trchb{{\tr \,\chib}}

\def\atrch{\atr\chi}
\def\atrchb{\atr\chib}

\def\rhod{\,\dual\hspace{-2pt}\rho}

\def\Ric{\mbox{ \bf Ric}}
\def\fb{\protect\underline{f}}
\def\err{{\mbox{Err}}}
\def\ov{\overline}

\def\f12{\frac 1 2}
\def\lab{\label}

\def\bsplit{\begin{split}}

\newcommand{\Mext}{{\,{}^{(ext)}\mathcal{M}}}

\newcommand{\Mint}{{\protect \,{}^{(int)}\mathcal{M}}}

\def\Mint{{\protect \, ^{(int)}\MM}}

\def\Mtop{{ \, ^{(top)} \MM}}
\def\Sitop{{\,^{(top)} \Si}}

\def\qf{\frak{q}}
\def\qfb{\protect\underline{\qf}}

\def\Jk{\mathfrak{J}}

\DeclareFontFamily{U}{mathx}{\hyphenchar\font45}
\DeclareFontShape{U}{mathx}{m}{n}{
      <5> <6> <7> <8> <9> <10>
      <10.95> <12> <14.4> <17.28> <20.74> <24.88>
      mathx10
      }{}
\DeclareSymbolFont{mathx}{U}{mathx}{m}{n}
\DeclareFontSubstitution{U}{mathx}{m}{n}
\DeclareMathAccent{\widecheck}{0}{mathx}{"71}

\def\Hbc{\widecheck{\Hb}}

\def\Gac{\widecheck{\Ga}}
\def\Rc{\widecheck R}

\def\squared{\dot{\square}}

\DeclareFontFamily{U}{mathx}{\hyphenchar\font45}
\DeclareFontShape{U}{mathx}{m}{n}{
      <5> <6> <7> <8> <9> <10>
      <10.95> <12> <14.4> <17.28> <20.74> <24.88>
      mathx10
      }{}
\DeclareSymbolFont{mathx}{U}{mathx}{m}{n}
\DeclareFontSubstitution{U}{mathx}{m}{n}
\DeclareMathAccent{\widecheck}{0}{mathx}{"71}

\def\und{\underline}

\def\uext{^{(ext)} u}
\def\uint{^{(int)} u}

\def\Ric{\mbox{Ric}}

\begin{document}

\title{Brief introduction to the nonlinear stability of  Kerr}
\author{Sergiu Klainerman, J\'er\'emie Szeftel}

\maketitle

\begin{abstract}
This a  brief introduction to   the  sequence of works  \cite{KS:Kerr}, \cite{GKS-2022},  \cite{KS-GCM1}, \cite{KS-GCM2}  and 
\cite{Shen}  which establish the nonlinear stability  of Kerr black  holes with  small angular momentum. We are  delighted 
to dedicate this article  to  Demetrios Christodoulou for whom  we both  have great admiration.
 The first author  would also like  to thank Demetrios  for the magic  moments  of friendship, discussions   and collaboration  he enjoyed 
  together  with him.
\end{abstract}


\section{Kerr stability conjecture}



\subsection{Kerr spacetime}
\lab{section:Kerr}


Let  $(\KK(a, m), \g_{a,m})$  denote   the family of  Kerr spacetimes
   depending  on the  parameters $m$  (mass)  and $a$ (with $J=am$    angular momentum).
 In Boyer-Lindquist  coordinates  the Kerr metric is given by
   \begin{equation}\label{zq1}
\g_{a,m}=-\frac{q^2\Delta}{\Sigma^2}(dt)^2+\frac{\Sigma^2(\sin\theta)^2}{|q|^2}\left(d\phi-\frac{2amr}{\Sigma^2}dt\right)^2 +\frac{|q|^2}{\Delta}(dr)^2+|q|^2(d\theta)^2,
\end{equation}
where
\begin{equation}\label{zq2}
\begin{cases}
&\Delta=r^2+a^2-2mr, \qquad  q=r+i a \cos \th,\\
&\Sigma^2=(r^2+a^2)|q|^2+2mra^2(\sin\theta)^2=(r^2+a^2)^2-a^2(\sin\theta)^2\Delta.
\end{cases}
\end{equation}
The asymptotically flat\footnote{That is they approach the Minkowski metric for large $r$.}  metrics    $\g_{a,m}$ verify   the Einstein vacuum equations
\bea
  \lab{eq:EVE-intro}
  \Ric(\g)=0,
  \eea  
     are   stationary and axially symmetric\footnote{That is  $\KK(a,m)$ possess two Killing vectorfields: the stationary  vectorfield    $\T=\pr_t$,  which is time-like   in the 
  asymptotic region, away  from the horizon, and  the  axial symmetric Killing field $\Z=\pr_{\phi}$.},  possess  well-defined  event horizons $r=r_+$ (the largest  root  of $\De(r)=0$), domain of outer communication $r>r_+$ and  smooth  future  null infinity $\II^+$  where  $r=+\infty$.   The metric can be extended  smoothly inside  the black hole region, see  Figure \ref{fig:penrosediagramofkerr}.   The boundary  $r=r_-$  (the smallest root  of $\De(r)=0$) inside the black hole region is a Cauchy horizon across which predictability fails\footnote{Infinitely many smooth extensions are possible beyond the  boundary.}. 
  
 \begin{figure}[h!]
\centering
   \includegraphics[width=5in]{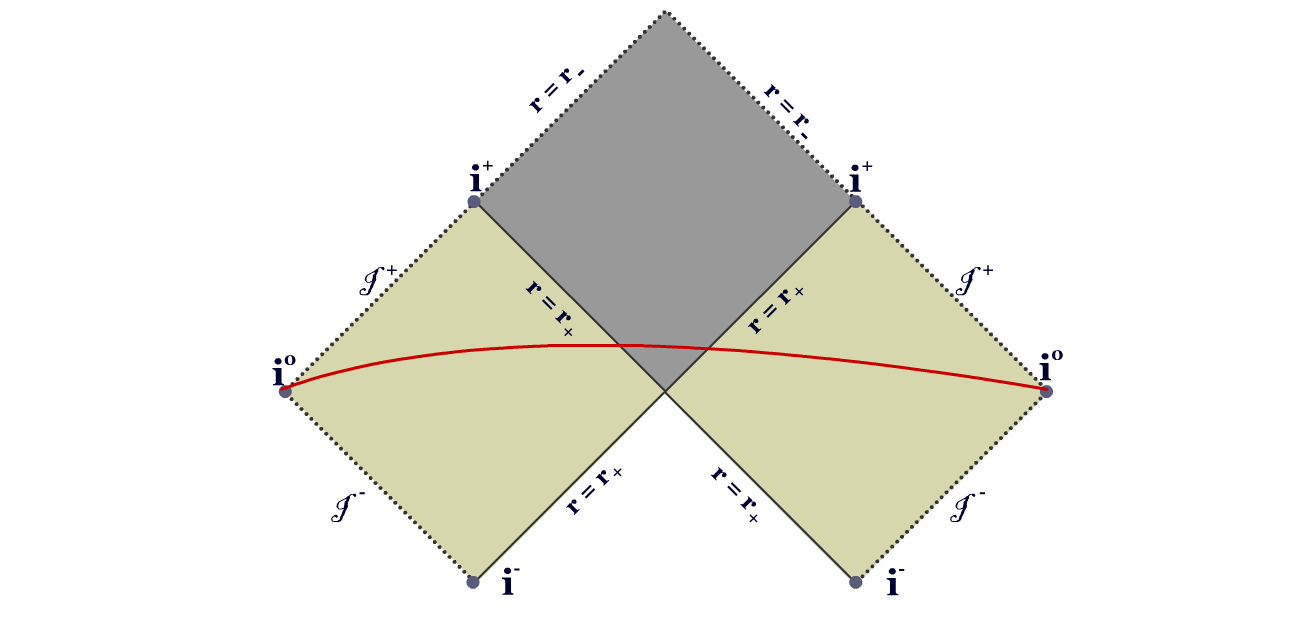}  
\caption{\footnotesize{Penrose diagram of Kerr for $0<|a|<m$. The surface  $r=r_+$, the larger root  of $\De=0$,  is the event horizon of the black hole, $r>r_+$ the domain of outer communication,   $\II^+$  is  the future null infinity, corresponding to $r=+\infty$. }}\lab{fig:penrosediagramofkerr}
\end{figure} 

Here are some of the most important  properties of $\KK(a,m)$:
\begin{itemize}
\item   $\KK(a, m)$  possesses  a canonical  family of  null  pairs, called \textit{principal null pairs},   of the  form
    $(\la e_4, \la^{-1} e_3) $, with $\la>0$ an arbitrary scalar function, and 
 \bea
 \lab{eq:PrincipalNull}
\bsplit
e_4&=\frac{r^2+a^2}{|q|^2} \pr_t +\frac{\De}{|q|^2} \pr_r +\frac{a}{|q|^2} \pr_\phi, \qquad
e_3  =\frac{r^2+a^2}{\De} \pr_t -\pr_r +\frac{a}{\De} \pr_\phi.
\end{split}
\eea 
\item  The horizontal structure, perpendicular to $e_3, e_4$, denoted $\HH$,  is spanned by the  vectors 
\bea
\lab{eq:canonicalHorizBasisKerr-intro}
e_1=\frac{1}{|q|}\pr_\th,\qquad e_2=\frac{a\sin\th}{|q|}\pr_t+\frac{1}{|q|\sin\th}\pr_\phi.
\eea
The distribution generated by $\HH$ is non-integrable  for $a\neq 0$.
\item The  horizontal  structure  $( e_3,  e_4, \HH)$  has   the remarkable property  that all components of the Riemann curvature  tensor  $\R$, decomposed  relative  to them,    vanish  with the exception of those which  can be deduced  from
$
\R(e_a, e_3, e_b, e_4).
$

 \item $\KK(a, m)$  possesses  the Killing vectorfields $\T, \Z$ which,  in BL coordinates, are given by  
 $\T=\pr_t, \Z=\pr_\phi$. 
 
  \item In addition to the symmetries  generated  by $\T, \Z$, 
   $\KK(a, m)$  possesses   also  a  non-trivial  Killing tensor\footnote{Given by  the expression $ \C=-a^2 \cos^2\th \g+O$,  $O=|q|^2 \big( e_1\otimes e_1+ e_2\otimes e_2\big)$.}, i.e. a symmetric   2-tensor $\C_{\a\b}$
   verifying the property
   $\D_{(\ga} \C_{\a\b)}=0$. The tensor carries the name  of its discoverer B. Carter, see \cite{Carter}, who  made use of it   to show that the   geodesic flow   in Kerr is integrable.  Its presence, in addition to $\T$ and $\Z$,  as a higher order symmetry,      is at the heart of what   Chandrasekhar, see \cite{Chand3},  called the most striking feature  of Kerr,   ``the separability of all the standard equations of mathematical physics in Kerr geometry''.

   \item  The Carter  tensor can be used to define the Carter   operator
   \bea
   \CC=\D_\a( \C^{\a\b}\D_\b),
   \eea
    a second order operator  which commutes  with $\square_{a,m}$. This property  plays a  crucial  role in  the proof of   our  stability  result,   Theorem \ref{MainThm-firstversion}, more precisely in  Part II of  \cite{GKS-2022}.
\end{itemize}


\subsection{Kerr stability conjecture}


The discovery of black holes, first  as  explicit solutions of  EVE and later as possible explanations of astrophysical phenomena\footnote{According to Chandrasekhar ``Black holes are macroscopic objects with masses varying from a few solar masses to millions of solar masses. To the extent that they may be considered as stationary and isolated, to that extent, they are all, every single one of them, described exactly by the Kerr solution. This is the only instance we have of an exact description of a macroscopic object. Macroscopic objects, as we see them around us, are governed by a variety of forces, derived from a variety of approximations to a variety of physical theories. In contrast, the only elements in the construction of black holes are our basic concepts of space and time. They are, thus, almost by definition, the most perfect macroscopic objects there are in the universe. And since  the general theory of relativity provides a single two parameter family of solutions for their description, they are the simplest as well.''}, 
 has not only revolutionized our understanding of the universe,   it  also gave mathematicians a monumental task: to test the physical reality of these solutions. This may seem nonsensical since physics tests the reality of its objects by experiments and observations and, as such, needs mathematics to formulate the theory and make quantitative predictions, not to test it. The problem, in this case, is that black holes are by definition non-observable and thus no direct experiments are possible.  Astrophysicists ascertain the presence of such objects through indirect observations\footnote{The physical reality of these objects was recently put to test  by  LIGO which is supposed to have detected the gravitational waves generated in the final stage of the in-spiraling of two black holes.   Rainer Weiss,    Barry C. Barish and  Kip  S. Thorne received the 2017 Nobel prize   for their ``decisive contributions'' in this respect.  The 2020 Nobel  prize in Physics  was  awarded  to   R. Genzel and A. Ghez    for providing  observational   evidence for the presence of super massive black holes in the center of our galaxy, and  to R.  Penrose  for his theoretical    foundational  work:   his  concept of a trapped surface and the proof of  his famous singularity theorem.}  and numerical experiments, but  both are limited in scope to the range of possible observations or 
 the  specific initial conditions  in which numerical  simulations are conducted.  One can rigorously check that the Kerr solutions have vanishing Ricci curvature, that is, their mathematical reality is undeniable.   But to be real in a physical sense, they have to satisfy certain properties which, as it turns out,  can  be neatly formulated in unambiguous mathematical language.
 Chief among them\footnote{Other  such properties concern the  rigidity of the Kerr family, see \cite{IK:rigidity} for a  current survey,   or  the dynamical formations  of black holes from regular configurations,  see  the \cite{Chr-BH},  \cite{KlRod}  and the introduction to \cite{An}  for an  up to date account of more recent results.} is the problem of stability, that is, to show that if the precise initial data corresponding to Kerr are perturbed a bit, the basic features of the corresponding solutions do not change much\footnote{If  the Kerr family would be  unstable  under  perturbations, black holes  would be nothing more than mathematical artifacts.}.

   {\bf Conjecture} (Stability of Kerr conjecture).\,\,{\it  Vacuum, asymptotically flat,  initial data sets, 
   sufficiently close to  $Kerr(a,m)$, $|a|/m <1$,  initial data,  have  maximal developments with complete
   future null infinity  and with   domain of outer communication\footnote{This presupposes the existence of an event horizon. Note that the  existence of such an event horizon   can only be established  a posteriori,  upon the completion of the proof of the conjecture.}   which
   approaches  (globally)  a nearby Kerr solution.}   

\bigskip


\subsection{Resolution of the conjecture for slowly rotating black holes} 
\lab{section:MainResult}



\subsubsection{Statement of the main result}   


The goal of this article is to  give a short introduction   to  our   recent    result in which we   settle  the  conjecture  in the case   of  slowly rotating  Kerr black holes. 
   \begin{theorem}
\lab{MainThm-firstversion}
The future globally hyperbolic   development  of  a general,   asymptotically  flat,    initial data set, sufficiently close  (in a suitable  topology)  to a   $Kerr(a_0, m_0) $   initial data set,  
  for sufficiently  small $|a_0|/m_0$,  has a complete    future null infinity  $\II^+$ and converges 
in  its causal past  $\JJ^{-1}(\II^{+})$  to another  nearby Kerr spacetime $Kerr(a_f, m_f)$ with parameters    $(a_f, m_f)$ close to the initial ones $(a_0, m_0)$.
 \end{theorem} 
 
 \begin{figure}[h!]
 \lab{fig0-introd}
\centering
\includegraphics[scale=0.6]{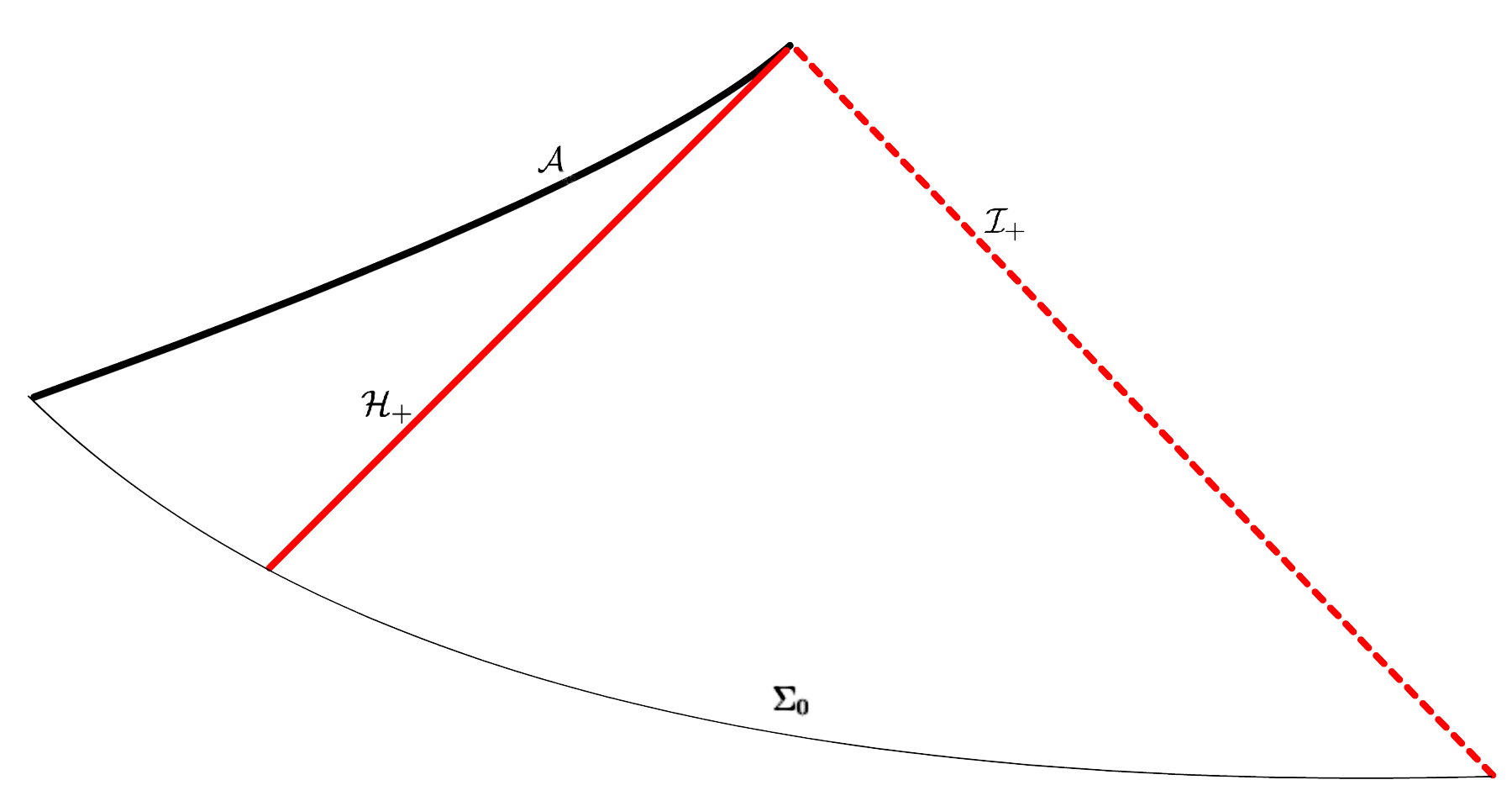}
\caption{\footnotesize{The Penrose diagram of the  final   space-time in  Theorem \ref{MainThm-firstversion} with initial hypersurface $\Si_0$,  future space-like boundary $\AA$, and $\II^+$ the complete future null infinity.  The hypersurface $\HH_+$} is the future event horizon of  the final Kerr.}
\end{figure} 

 The precise version of   the result,    all the main  features of the   architecture of its proof,  as well    as detailed proofs  for  most of the main steps  are to be found in  \cite{KS:Kerr}.  
  The full proof relies  also on  our joint work  \cite{GKS-2022} with E. Giorgi, our papers \cite{KS-GCM1}, \cite{KS-GCM2} on GCM spheres, and the extension \cite{Shen} to GCM hypersurfaces by D. Shen.


\subsubsection{Brief comments on  the  proof}   
\lab{section:Brief}


 We  will discuss the main ideas of the proof in more  details  in section \ref{section:Proof}.  It pays  however  to   give already a     graphic  sense   of   the main building  blocks  of  our  approach,  which we call  \textit{general covariant modulated}(GCM), admissible  spacetimes.

 \begin{figure}[h!]
 \centering
\includegraphics[scale=0.9]{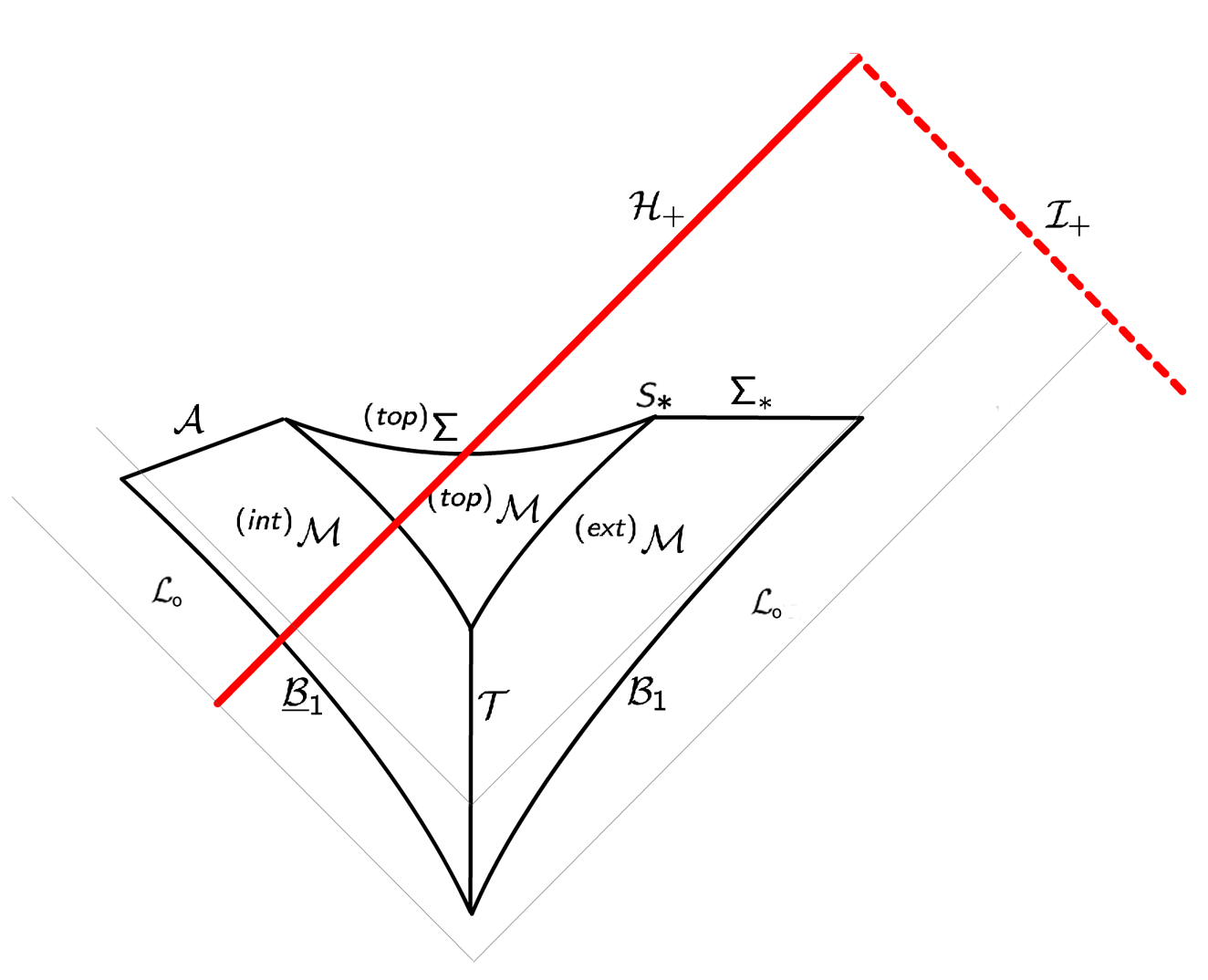}
\caption{\footnotesize{The Penrose diagram of  a finite GCM admissible    space-time $\MM=\Mext\cup\Mtop\cup\Mint$. The future  boundary 
 $\Si_*$       initiates at the  GCM sphere $S_*$.      The  past boundary of $\MM$,  $\BB_1\cup\BBb_1$, is  included in  the  initial layer $\LL_0$, in which the spacetime is assumed  given.   }}
\lab{fig1-introd1} 
\end{figure}

 The   main  features of these  spacetimes $\MM=\Mext\cup\Mtop\cup\Mint$  are   as follows:

\begin{itemize}
\item  The capstone   of the entire construction is the  sphere $S_*$,  on the  future boundary  $\Si_*$ of $\Mext$,  which verifies a  set of  specific extrinsic and intrinsic  conditions denoted by the acronym GCM. 

\item The  spacelike hypersurface $\Si_*$,   initialized at $S_*$,    verifies a set of   additional GCM  conditions. 

\item  Once $\Si_*$ is specified the whole  GCM admissible spacetime $\MM$   is  determined by  a more  conventional  construction, based on  geometric transport type equations\footnote{More precisely $\Mext$  can be determined  from $\Si_*$ by  a specified  outgoing  foliation terminating in the timelike boundary $\TT$,  $\Mint$ is determined from $\TT$ by a  specified  incoming 
 one, and $\Mtop$  is a complement of $\Mext\cup\Mint$   which makes $\MM$ a  causal domain.}.  
 
 \item The construction, which   also allows us to specify  adapted null frames\footnote{In our work we  prefer to   talk about horizontal structures, see the brief  discussion  in section \ref{sec:nonintegrabilityhorizontalstructure}. Another important  novelty in the proof of Theorem \ref{MainThm-firstversion}  is that it relies  on non-integrable horizontal structures, see section \ref{sec:nonintegrabilityhorizontalstructure}.},  is  made possible by  the   covariance properties  of the Einstein  vacuum equations.  
 
 \item  The past boundary $\BB_1\cup\BBb_1$ of $\MM$, which is itself  to be constructed,  is included in the initial layer $\LL_0$ in which the  spacetime is assumed 
 to be known\footnote{The passage form  the initial data specified  on $\Si_0$ to  the initial layer  spacetime $\LL_0$,   can be  justified 
  by arguments similar to those of  \cite{KlNi1} \cite{KlNi2},   based  on the  methods introduced in \cite{Ch-Kl}.}, i.e. a small vacuum  perturbation of a Kerr solution.  
\end{itemize}

The proof   of Theorem  \ref{MainThm-firstversion}  is  centered   around    a  limiting  argument for  a continuous family  of such spacetimes  $\MM$  together  with   a set   of bootstrap assumptions (BA)  for the connection and curvature coefficients,  relative to  the  adapted frames.   Assuming that  a given  finite, GCM admissible, spacetime   $\MM$   saturates  BA  we  reach a contradiction as follows:
\begin{itemize}
\item First  improve BA   for   some of the components of the curvature tensor with respect to the frame. These  verify  equations (called Teukolsky equations)   which decouple, up to terms quadratic in the perturbation, and    are treated by wave equations  methods.
\item Use the information provided by these curvature coefficients  together with the  gauge choice  on $\MM$,  induced  by the GCM   condition on $\Si_*$,  to    improve BA for all other  Ricci and curvature components.
\item Use  these improved estimates to extend $\MM$ to a strictly  larger   spacetime $\MM'$ and then construct a new GCM  sphere  $S_*'$,  a new    boundary   $\Si_*'$ which initiates on $S_*'$,   and a new GCM admissible spacetime   $\MM'$, with $\Si_*' $ as  boundary,  strictly larger than $\MM$.
\end{itemize}

\begin{remark}
The  critical new feature  of this argument    is the  fact that  the new  GCM sphere  $S_*'$  has to be constructed as a co-dimension 2  sphere in $\MM'$
with no  reference to the initial conditions\footnote{See   a more  detailed discussion in section \ref{sec:GCMadmissiblespacetimes}.}. This   construction   appears  first in \cite{KS} in a polarized situation. The general construction   appears in the GCM   papers  \cite{KS-GCM1},    \cite{KS-GCM2}.
The construction of $\Si_*'$ from  $S_*'$  appears first in \cite{KS} in the polarized case. The general construction  used in our work  is    due to D. Shen \cite{Shen}.
\end{remark}

  
 \section{Linear and nonlinear  stability}


  
\subsection{Notions of nonlinear  stability} 
\lab{section:NotionsofLNL-stability} 


Consider  a stationary solution $\phi_0$ of a  nonlinear evolution equation 
\bea
\lab{eq:Nonl-NN}
\NN[\phi]=0.
\eea 
  There are two  distinct notions of stability, \textit{orbital stability}, according to which  small perturbations    of $\phi_0$    lead to solutions $\phi$ which remain close for all time, and  \textit{asymptotic stability (AS)} according to which the perturbed solutions converge as $t\to \infty$ to $\phi_0$. In the case where $\phi_0$ is non trivial, there is a third notion, which we call \textit{asymptotic orbital stability (AOS)}, to describe the fact that  the perturbed solutions may converge to a different stationary solution. This happens if $\phi_0$ belongs to a multi parameter smooth family of stationary solutions, or by applying a gauge transform to $\phi_0$ which keeps the equation invariant\footnote{In the case of Kerr, both cases are present as we shall see below.}. 
  
  For quasilinear equations\footnote{Orbital stability can be  established directly (i.e. without   establishing the stronger version) only  in rare occasions,   such as for  hamiltonian equations with weak  nonlinearities.}, such as  EVE,   a proof of stability means necessarily  AS  or AOS   stability.      Both    require    a    detailed  understanding of the \textit{ decay properties}    of the   linearized equation, i.e. 
   \bea
   \lab{eq:Linearized-LL}
   \LL[\phi_0]\psi =0,
   \eea
    with $\LL[\phi_0]$ the Fr\'echet  derivative   $\NN'[\phi_0]$. This   is,  essentially, 
  a linear hyperbolic  system  with variable coefficients  which,  typically,  presents   instabilities\footnote{In  unstable situations    \eqref{eq:Linearized-LL} may have  
  exponentially growing solutions, see for example \cite{DKSW}.}.

  In the exceptional  situation,  when  stability can ultimately  be established,    one  can  tie  all   the  instability   modes to  the    following   properties of the nonlinear equation:  
   \begin{enumerate}
\item[M1.]  If  $\phi_\la $ is a family of stationary  solutions,  near $\phi_0$,     verifying    $\NN[\phi_\la ]=0$.  Then  $\psi_0 =(\frac{d}{d\la} \phi_\la)_{\la=0} $ verifies 
  $\NN'[\phi_0] \psi_0=0$, i.e. $\psi_0$ is a nontrivial, stationary,  bound state  of the linearized  equations \eqref{eq:Linearized-LL}.  
   
\item[M2.] If   $\Phi_\la $ is a smooth  family of  diffeomorphisms   of the background manifold,  
$\Phi_0=I$,   such that $ \NN[\Phi^*_\la(\phi_0)]=0$. Then
$\Psi_0=\big(  \frac{d}{d\la}  \Phi^*_\la(\phi_0)\big)_{\la=0}$ verifies $\NN'[\phi_0]\Psi_0=0$, i.e. $\Psi_0$ is also
 a stationary bound state  of the linearized equation \eqref{eq:Linearized-LL}.  
\end{enumerate}
  
  These  linear  instabilities  are responsible for the fact   that a small perturbation of  the   fixed stationary solution $\phi_0$  may not converge to  $\phi_0$ but to  another nearby  stationary solution\footnote{The methodology of tracking this asymptotic final state, in general different from $\phi_0$, is usually referred to as modulation.  See   for example  \cite{Ma-Me},\cite{Me-R}  for how modulation theory  can be used   to deal with some examples of scalar nonlinear dispersive equations.}.
  
  To prove the asymptotic  convergence   of  $\phi$   to  a  final state $\phi_f$, different  form $\phi_0$,  we  need  to establish  sufficiently strong rates  of decay\footnote{To control the nonlinear terms of the equation.} for $\phi-\phi_f$.          Rates of decay however are strongly  \textit{coordinate dependent},   i.e.  dependent on   the  choice of  the   diffeomorphism (or gauge) $\Phi$ in which decay is measured.  
   Thus, to prove  a  nonlinear stability,   result we need to know both  the final   state $\phi_f$ and the coordinate system $\Phi_f$   in which sufficient decay,  and thus   convergence  to $\phi_f$, can be established. The difficulty  here is that neither $\phi_f$ nor $\Phi_f$   can be determined a-priori  (from  the  initial perturbation),  they have  to emerge dynamically  in the process of convergence. Moreover,  in all   examples  involving  nonlinear wave equations in $1+3$ dimensions, the nonlinear terms  have to also cooperate, that is  an appropriate  version  of the so called  null condition  has to be verified.

To  summarize,  given a nonlinear system  $\NN[\phi]=0$ which possesses both a smooth family of stationary solutions $\phi_\la$   and   a smooth family of invariant diffeomorphisms a proof of the nonlinear stability of $\phi_0$  requires the following ingredients:
\begin{itemize}
\item  The only non-decaying modes of the linearized  equation $\LL(\phi_0)\psi=0$ are those due to  the items M1--M2 above. In particular there are no exponentially growing  modes.

\item A  dynamical construction of both the final state $\phi_f$ and the final gauge $\Phi_f$ in which convergence to the final state takes place.

\item The nonlinear  terms in   the equation
\beaa
\LL[\phi_f] \psi = N(\psi)
\eeaa
obtained by expanding  the equation $\NN[\phi]$  near $\phi_f$, in the gauge  given by  the diffeomorphism $\Phi_f$, 
has  to verify  an appropriate version of the null condition.
\end{itemize}


\subsection{The case of  the Kerr family}


 The issue of the stability of the Kerr    family 
      has been at the center of attention  of GR physics and  mathematical relativity for   more than half a century, ever since their  discovery   by Kerr in \cite{Kerr}.   In this case  we have to deal 
       not only with  a $2$-parameter  family of  solutions,  corresponding to   the parameters $ (a, m)$, but also  with the entire group of diffeomorphisms\footnote{Indeed, according to the covariant properties of the Einstein vacuum equations  we  cannot distinguish between $\g$ and $\Phi^*\g$,   for any diffeomorphism $\Phi$ of $\MM$.} of  $\MM$.    In what follows we  try to discuss  the main  difficulties of the problem. In doing that  it helps to  compare these  to those  arising    in the simplest  case  when $a=m=0$, i.e.  stability of Minkowski.

       
\subsection{Stability of Minkowski space} 
\lab{section:stabilityMink}


   Until very recently  the only   space-time for which  full nonlinear  stability had been established was the Minkowski space, see  \cite{Ch-Kl}. The proof is based  on  some     important  PDE  advances    of late last century:
\begin{enumerate}
\item[(i)]  Robust    approach, based on the  vectorfield method,  to   derive quantitative   decay  based on generalized energy estimates and commutation with (approximate)  Killing and conformal Killing  vectorfields. 

\item[(ii)]   The  \textit{null condition} identifying the deep  mechanism for     nonlinear stability, 
i.e.  the specific structure of the nonlinear terms   which   enables stability  despite the low
  decay of the  perturbations. 
  
  \item[(iii)]    Elaborate bootstrap argument according to which  one makes  educated assumptions 
about the behavior of  solutions  to  nonlinear wave equations  and then proceeds, by a long sequence of a-priori estimates,  to       show that   they are in fact 
 satisfied. This  amounts to a \textit{conceptual linearization}, i.e. a method  by which the equations become, essentially, linear\footnote{Note that in the context of  EVE, and other quasilinear hyperbolic systems, this    differs substantially from the usual notion of linearization around a fixed background.}  without   actually linearizing them.
 \end{enumerate}
  The main innovation in the  proof in \cite{Ch-Kl}  is   the choice  of  an  appropriate   gauge   condition,  readjusted dynamically  through  the  convergence process,   by  a  continuity argument,
  which  allows one to separate   the curvature  estimates,   treated by hyperbolic methods,   from the    estimates for the connection coefficients. The key point is that   these   latter  verify  transport or elliptic equations in which the curvature terms appear as sources.  Thus both  the curvature  components and 
  connection coefficients  can be  controlled by a bootstrap argument.
    The  gauge condition is based on  the  constructive choice  of a maximal  time function $t$ and two outgoing optical functions  $\uint$\footnote{The interior optical function is initialized  on  a  timelike  geodesic  from the  initial hypersurface.} and $\uext$\footnote{The exterior  optical  function  $\uext$ is initialized  on   the last slice $t=t_*$, by  the  construction of a    foliation (inverse lapse foliation)  initialized at space-like infinity.  It   is thus readjusted  dynamically  as $t_*\to \infty$.}  covering  the interior and exterior  parts of  the spacetime.

  Another novelty of \cite{Ch-Kl}  was the reliance on    null frames  adapted to   the  $S$-foliations   induced    by the level surfaces of $t$ and $u$. These  define   integrable horizontal structures (in the language of   part I of  \cite{GKS-2022}),   by contrast with the non-integrable ones  used in the proof of Theorem      \ref{MainThm-firstversion} and discussed in section \ref{sec:nonintegrabilityhorizontalstructure}. The functions $ t, u$ and  this  integrable horizontal structure can be used to define approximate Killing vectorfields used in estimating the curvature.

   
\subsection{Main difficulties}
    \lab{subs:MainDiffic}
    
    
  There are a few   major obstacles in passing from the stability of Minkowski  to that of  Kerr:
   \begin{enumerate}
 \item The first one was     already discussed  in  section \ref{section:NotionsofLNL-stability} in the general context of the stability of a stationary solution $\phi_0$.  In the case  when $\phi_0$ is trivial    there are no nontrivial bound states    for the linearized problem and thus   we expect  that the final state   does actually coincide  with  $\phi_0$.  This  is precisely  the case  for   the special member  of the Kerr family 
  $a=m=0$, i.e. the Minkowski  space\footnote{Note however  that  even though the linearized system around Minkowski  does not contain instabilities,   the proof of the nonlinear stability of the Minkowski space  in \cite{Ch-Kl}    takes into account (in a fundamental way!)   general covariance. Indeed  the presence of  the   ADM mass affects the causal structure  of  the far, asymptotic,  region of the perturbed space-time.} $(\RRR^{1+3}, m)$.         On the other hand, in perturbations of Kerr,  general covariance affects the entire construction of the spacetime. In the proof of Theorem \ref{MainThm-firstversion}  the crucial  concept  of a GCM admissible  spacetime  is meant to  deal  with  both  finding the final parameters and the  gauge in which convergence to the final state takes place.

\item A  fundamental insight   in the stability of the Minkowski space  was that the Bianchi identities  decouple at first order from the null structure equations which allows one to control curvature first, as a Maxwell type system (see \cite{Ch-Kl0}), and then proceed with the rest of the solution. This  cannot work for  perturbations of Kerr due to the fact that some of the null components\footnote{With respect to  the  \textit{principal null directions of Kerr,} i.e   a distinguished  null  pair which   diagonalizes  the  full curvature tensor, the middle component $P=\rho+ i\rhod$ is nontrivial.} of the  curvature tensor  are non-trivial in Kerr.
 
 \item    Even if  one succeeds in tackling the above mentioned issues, there are still major obstacles  in understanding the  decay properties of the solution.  Indeed,  when one  considers  the simplest,    relevant,
  linear equation  on a fixed Kerr background,  i.e. the scalar wave equation $\square_{a,m} \psi=0$,  one 
  encounters serious difficulties  to  prove decay.   Below is a  very short description of these:
       \begin{itemize}
        \item \textit{The problem of   trapped null geodesics.}   This  concerns the existence of null geodesics\footnote{In the Schwarzschild case, these geodesics
are associated with the celebrated photon sphere $r=3m$.} neither crossing the event horizon nor escaping to 
null infinity, along which solutions   can concentrate for arbitrary long
times.  This   leads  to degenerate  energy-Morawetz estimates  which require a very delicate analysis.
       
\item \textit{The   trapping   properties of the     horizon.}       
    The horizon  itself  is ruled  by null  geodesics, which do not     communicate with  null infinity and can thus concentrate energy.         This problem was solved   by understanding      the  so called     red-shift effect associated to  the event horizon,  which     counteracts this type of  trapping.

\item \textit{The problem of superradiance.} This  is  the failure of the stationary Killing 
field $\T=\partial_t$
to be everywhere timelike in the domain of outer communication\footnote{The stationary  Killing  vectorfield $\T$ is timelike only outside of the so-called ergoregion.}, and thus, of the
associated  conserved energy to be positive. Note that this problem is absent in Schwarzschild and, in general,   for axially symmetric solutions of EVE. In both cases however there still is a degeneracy along the horizon.
\item   \textit{Superposition problem.}       This is  the  problem  of combining the estimates  in the near  region, close to the horizon, (including the  ergoregion  and trapping) with estimates in the asymptotic region, where the spacetime is close to   Minkowski. 
\end{itemize}

\item Though, as seen above,  the analysis of  the scalar wave equation in Kerr  presents formidable difficulties,  it is itself just  a vastly simplified model problem.  A more realistic  equation is the so called spin 2 wave equation, or Teukolsky equation,  which presents many new challenges\footnote{Unlike  the  scalar wave
 equation $\square_{a,m} \psi=0$, which is conservative, the Teukolsky equation is not, 
 and  we  thus  lack  the most basic ingredient  in  controlling     the   solutions of the equation, i.e. energy estimates.}.

   \item   The full linearized system, whatever  its  formulation,  presents   many additional difficulties   due to its  complex  tensorial  structure and  the  huge gauge covariance of the  equations\footnote{As mentioned  earlier, rates of decay are heavily dependent on a proper choice of gauge, thus affecting the issue of convergence.}.   The crucial  breakthrough in this regard is the  observation, due to Teukolsky \cite{Teuk},  that the  extreme components of the  linearized curvature tensor are both gauge invariant (see below in  section \ref{sect:LinearStab}) and  verify decoupled spin 2  equations, that is  the Teukolsky equations mentioned above.   

\item   A crucial simplification
 of the linear theory, by comparison to the  full nonlinear case, is that  one   can  separate the treatment of the gauge  invariant   extreme  curvature components form all the other gauge invariant quantities.
  In  the nonlinear case this separation is no longer true, all quantities need to be treated simultaneously.  Moreover,  methods based on  separation of variables, developed  to treat scalar and and spin 2 wave  equations in Kerr,    are   incompatible  with the nonlinear setting which requires, instead,  robust methods   to derive decay.  
\end{enumerate}

   
\subsection{Linear stability}
\lab{sect:LinearStab}


Linear stability  for the vacuum  equations is formulated in the following way.  Given
the Einstein tensor $\G_{\a\b}=\R_{\a\b}-\frac 1 2 \R \g_{\a\b}$  and a stationary solution 
$\g_0$, i.e. a fixed Kerr metric, one   has to solve the  system of equations
 \bea
\lab{eq:Lin. Grav}
\G'(\g_0) \, \de \g =0.
\eea
The  covariant properties of the  Einstein equations, i.e. the equivalence between a solution $\g$ 
 and $\Phi^*(\g)$, leads  us to identify
    $\de\g$  with   $\de\g+\Lie_X (\g_0)$ for arbitrary  vectorfields    $ X$    in $\MM$, i.e.      
   \bea
   \lab{eq:linearcovariance}
 \de\g\equiv  \de\g+\Lie_X (\g_0).
   \eea
   We can now attempt to  formulate a version of linear stability for \eqref{eq:Lin. Grav}, loosely related to   the   nonlinear stability of Kerr    conjecture,  as follows.
   \begin{definition}
  By   linear stability of the Kerr  metric $\g_0$   we understand a result  which achieves the following:
  
   Given  an  appropriate initial data  for  a perturbation  $\de\g$,   find a  vectorfield $X$ such that,
   after projecting  away the   bound states   generated by the parameters $a, m$
    according to  M1--M2  in section \ref{section:NotionsofLNL-stability},   all solutions  of the form  $\de\g+\Lie_X\g_0$      to \eqref{eq:Lin. Grav},   decay, relative  to  an  appropriate  null  frame of $\KK(a_0, m_0)$,    
      sufficiently fast  in time.
   \end{definition} 

   \begin{remark}
   The definition above is necessarily  vague.       What    is the meaning of sufficiently fast?       In fact    various  components of the metric $\de\g$, relative   to  the canonical null frame  of $\KK(a, m)$,   are expected to decay at different  polynomial rates, some of which  are  not even integrable.      Unlike in the nonlinear context, where one needs precise   rates of    decay  for components   of the curvature tensor  and Ricci coefficients,  as well as    their derivatives,   to be able to control the nonlinear terms, in linear  theory     any type of  nontrivial   control of  solutions  may be  regarded as satisfactory\footnote{This is in fact what happens in the literature  on linear theory. Thus, for example,  in their well known  linear stability result around Schwarzschild  \cite{D-H-R}, the authors  derive   satisfactory results (compatible    with  what is needed in nonlinear theory)  for  components of the  curvature  tensor, and some  Ricci coefficients, but not all.}.  
     Thus linear stability, as formulated above,   can only be regarded as 
    a vastly simplified  model problem. Nevertheless  the study of linear stability of the Kerr family  has turned    out to be    useful in various ways, as we shall see below. 
   \end{remark}

   Historically,     the following  versions of linear stability  have been   considered.
   \begin{itemize}
   \item[(a)]\textit{Metric Perturbations.} At the level of the metric  itself, i.e.  as above in \eqref{eq:Lin. Grav}.  
   \item[(b)]  \textit{Curvature Perturbations.}  Via  the Newman-Penrose (NP) formalism,  based  on null frames. 
   \end{itemize}
   The strategy  followed  in  both  cases\footnote{In the article we  refer mainly to  the curvature perturbation approach.}     is:
   \begin{itemize}
  \item  Find  components  of the metric (in  case (a)) or curvature tensor  (in case (b)), invariant  with respect to linearized gauge transformations \ref{eq:linearcovariance}  which verify   simple, decoupled, wave equations.   
   The main insight of this type was  the discovery, by  Teukolsky \cite{Teuk},  in the context of  (b) above,  that the extreme components of the  linearized curvature tensor  verify   both  these properties. 
  \item Analyze  these components   by showing one of the following:
  \begin{itemize}
  \item    There  are no exponentially  growing modes. This  is known as mode stability.
  \item Boundedness  for all time.
  \item      Decay (sufficiently fast)    in time.
  \end{itemize} 
  \item Find a linearized  gauge condition, i.e. a vectorfield $X$, such that  all remaining (gauge dependent) components (at the metric or curvature level)   inherit  the property mentioned above: no  exponentially growing modes,  boundedness, or decay in time.  In the physics literature  this  is  known as the problem of  reconstruction.  
   \end{itemize}

   
\subsubsection{Mode stability}  
\lab{section:ModeStab}


   All results on the  linear  stability of Kerr   in the   physics literature  during the  10-15 years after  
    Roy Kerr's  1963  discovery,   often called  the ``Golden Age of Black Hole Physics'',   are based 
   on mode decompositions. One makes use of the separability\footnote{See  discussion in  section \ref{section:Kerr}.} of the  linearized  equations,  more precisely the Teukolsky equations,   on a fixed Kerr background, to derive   simple  ODEs  for the corresponding modes.     One can then show, by ingenious methods,  that these modes  cannot exhibit exponential growth.       The  most complete  result  of    this type  is due to Bernard Whiting  \cite{Whit}.

   The obvious limitation of these results    are as follows:
   \begin{itemize}
   \item  They  are far from   even  establishing    the  boundedness of  general solutions  to  the Teukolsky equations, let alone  to establish  quantitative decay  for the  
   general solutions.
   \item  Results based on mode decompositions   depend strongly on the specific symmetries   of   Kerr
    which   cannot be adapted to  perturbations  of Kerr.
   \end{itemize}
   
     Robust  methods  to deal  with both issues have been developed  in the  mathematical community, based on the vectorfield method which we  discuss below.

 
\subsubsection{Classical vectorfield method} 


The  vectorfield method, as an analytic tool to derive   decay,   was first   developed in connection   with the wave  equation  in Minkowski space.    
As well known,  solutions of the wave equation  $ \square \psi=0$ in the Minkowski space  $\RRR^{n+1}$  
  both conserve energy  and decay uniformly in time   like $t^{-\frac{n-1}{2}}$. While conservation of energy  can be established  by a simple integration by parts,  and is thus robust to perturbations of the Minkowski metric,  decay   was first derived   either using  the   Kirchhoff formula  or by Fourier methods, which are manifestly not robust.  An integrated version of local energy  decay, based on an  inspired  integration by parts argument,   was first derived by  C. Morawetz \cite{Mor1}, \cite{Mor2}.  The first    derivation of decay 
   based  on   the  commutations properties of  $\square$ with     Killing and  conformal Killing  vectorfields    of Minkowski space    together with    energy conservation 
   appear in \cite{Kl-vectorfield} and \cite{Kl-vect}.   The same 
    method also  provides  precise  information about  the    decay properties  of  derivatives of solutions  with respect to   the standard null frame 
    of Minkowski space,   an important motivating factor in the discovery of the  null condition \cite{Kl-ICM}, \cite{Chr} and \cite{Kl-null}.    
    
      The crucial feature  of the  methodology  initiated by these papers,    to which we  refer      as the classical vectorfield method,    is that  it can be easily   adapted to perturbations of  the Minkowski space.  
      As such the method     has had numerous applications to nonlinear wave equations and played an important  role in  the proof of the nonlinear stability of Minkowski space, as discussed in section \ref{section:stabilityMink}. It  has also been applied    to   later versions  of the stability  of Minkowski   in  \cite{KlNi1}-\cite{KlNi2},   \cite{Lind-Rodn},  \cite{Bi}, \cite{Lind}, \cite{huneau}, \cite{HVas2}, \cite{graf}, and extensions of it   to  Einstein equation coupled with  various matter  fields in \cite{BiZi},  \cite{FJS}, \cite{Lind-T}, \cite{BFJT}, \cite{Wa},   \cite{Lf-Ma},   \cite{I-Pau}.

  
\subsubsection{New vectorfield method}
\lab{section:NewVf}
  
  
  To  derive    decay estimates    for solutions of wave equations  on a Kerr background one has to  substantially refine the classical vectorfield  method.
The  new vectorfield method  is an extension of the classical  method  which  compensates  for   the lack of  enough Killing and conformal Killing vectorfields   in Kerr   by   introducing,   new,  cleverly  designed,   vectorfields  whose deformation tensors have    coercive properties    in different regions of spacetime, not necessarily causal.   
The method  has    emerged  in the last 20  years   in connection to the study of boundedness and decay  for the    scalar wave equation in Schwarzschild and Kerr, see section \ref{section:QuantitativeDecay}  for more details.

   
\subsection{Model problems}  
\lab{sect:ModelProblems}


To  solve  the stability of Kerr conjecture  one has to   deal simultaneously  with all the difficulties  mentioned above.  This is, of course, beyond the abilities  of  mere  humans.  Instead  the problem was tackled in  a  sequence  of  steps based on  a variety  of simplified  model problems,  in  increasing  order of   difficulty.   To start with 
we can classify  model problems based on the following criteria:
\begin{enumerate}
\item Whether  the result refers to Schwarzschild i.e. $a=0$,  slowly rotating Kerr i.e. $|a|\ll m$   or   full non-extremal  Kerr $|a|<m$.
\item Whether the result  refers to  linear or nonlinear stability.
\item Whether the result, in linear theory,   refers  to scalar wave equation,  i.e. spin $0$,  Teukolsky equation,  i.e.  spin  $2$,  or the full linearized system.
\item Whether  the  stability   result, in linear theory,  is a mode stability result, a boundedness result  or one that establishes  quantitative decay. 
\end{enumerate}

      
\section{Short survey  of  model problems}

  
   We   give  below a short 
      outline of the main developments  concerning linear and nonlinear model problems for the Kerr stability problem, paying   special attention to  those which   had a measurable    influence on our work.

      
\subsection{Mode stability results} 
\lab{section:Modestability}
       
       
These are  mode stability  results, using the method of separation of variables,   obtained in the physics community  roughly during the period 1963-1990.  They  rely on what   Chandrasekhar 
        called     the most striking feature of Kerr  i.e.  ``\textit{the separability  of all the standard equations of mathematical physics in Kerr geometry''}. 
    
     \begin{enumerate}
 \item \textit{Regge-Wheeler (1957).}   Even before the discovery of the  Kerr solution physicists were interested in the mode  stability of  Schwarzschild space, i.e. $\KK(0, m)$. The first important    result  goes back to T. Regge and J.A Wheeler
       \cite{Re-W},  in  which they analyzed   linear, metric perturbations, of the Schwarzschild metric. They   showed that in a suitable gauge, equation \eqref{eq:Lin. Grav}  decouples into even-parity and
odd-parity perturbations, corresponding to axial and polar perturbations.    The  most important  discovery in that paper  is that of  the master Regge-Wheeler equation, a wave  equation  with a  favorable  potential, verified by an  invariant  scalar component $\phi$ of  the metric, i.e.
\bea
\lab{eq:Regge-W}
\square_{m}\phi= V\phi, \qquad  V=\frac{4}{r^2}\left(1-\frac{2m}{r}\right).
\eea
 where $\square_m$ denotes the  wave operator of the Schwarzschild   metric of mass $m$.  
 The R-W   study was completed by Vishveshwara \cite{Vishev} and Zerilli  \cite{Ze}.  A gauge-invariant formulation of \textit{metric perturbations}  was then given by Moncrief  \cite{Moncr}.

     \item  \textit{Teukolsky (1973).} The curvature perturbation approach, near Schwarzschild,    based on the  Newman-Penrose  (NP) formalism       was  first  undertaken  by Bardeen-Press \cite{Bar-Press}. This approach  was later extended to the Kerr family by Teukolsky  \cite{Teuk}, see also  \cite{P-T},   who made the  important  discovery that  the extreme  curvature components, relative   to  a principal null frame,  are  gauge invariant  and  satisfy  decoupled, separable,  wave equations. The equations, bearing the name of   Teukolsky,  are roughly   of the  form
     \bea
     \lab{Teuk}
     \square_{a,m} \psi=\LL[\psi]
     \eea
     where $\LL[\psi]$ is a first order linear operator   in $\psi$.
     
     \item \textit{Chandrasekhar (1975)}.  In \cite{Chand2}   Chandrasekhar   initiated  a  transformation theory relating  the  two   approaches.  He exhibited a transformation which connects  the Teukolsky  equations     to a       Regge-Wheeler  type equation.  In the particular case of Schwarzschild 
 the transformation  takes the Teukolsky equation       to the Regge-Wheeler  equation in \eqref{eq:Regge-W}.
     The   Chandrasekhar  transformation  was further  elucidated and extended  by R. Wald \cite{Wald} and recently by  Aksteiner and al \cite{Akst}.
     
      Though originally  it was    meant only   to    unify the   Regge-Wheeler approach with that of Teukolsky, the  Chandrasekhar transformation, and various   extensions  of it,  turn out  to  play an important   role in  the field. 
      
     \item \textit{Whiting (1989)}.   As mentioned before the full  mode stability, i.e. lack of exponentially growing modes,  for the Teukolsky equation on Kerr is  due  to  Whiting \cite{Whit} (see also \cite{Yacov},  in the case of the scalar wave equation, and \cite{AMPW} \cite{Teix} for stronger quantitive versions).  
       
       \item \textit{Reconstruction}. Once  we know  that the Teukolsky variables, i.e. the extreme components of the curvature tensor   verify mode stability, i.e. there are no exponentially growing modes, it still remains  to  deal with the problem of reconstruction, i.e. to  find a gauge relative to which all other components of the curvature 
     and Ricci coefficients enjoy the same property.  We refer the  reader  to Wald \cite{Wald}  and the   references within  for a treatment of this  issue in the physics literature.
 \end{enumerate}

     
\subsection{Quantitative decay for the scalar wave equation}
\lab{section:QuantitativeDecay}

     
As mentioned in section \ref{section:ModeStab},  mode stability     is far from establishing  even the  boundedness of solutions. To achieve that\footnote{The first   realistic boundedness result for solutions of the scalar wave equation in Schwarzschild appears in 
  \cite{K-Wald} based on a clever use of the energy method  which takes into account  the degeneracy of $\T$ at the horizon.}  and, more  importantly, to   derive realistic decay estimates, 
one needs  an entirely different approach based  on what we called  the  ``new vectorfield  method'' in  section \ref{section:NewVf}.    The method  has    emerged    in connection to the study of boundedness and decay  for the    scalar wave equation in $\KK(a,m)$,
 \bea
 \lab{eq:WaveKerr}
  \square_{a,m} \phi=0. 
  \eea

  The starting and most  demanding   part of the new method,  which appeared first in \cite{B-S1}, is  the derivation of  a global,  combined,  \textit{Energy-Morawetz}   estimate  which degenerates   in the trapping region.  Once an Energy-Morawetz estimate is established  one  can  commute  with the   Killing    vectorfields   of the background manifold, and  the so called red shift  vectorfield  introduced in  \cite{DaRo1},   to derive uniform bounds  for  solutions.  The most efficient  way  to also  get decay, and  solve the\textit{ superposition  problem} (see section \ref{subs:MainDiffic}), originating in \cite{Da-Ro3},  is  based on   the presence of   family of \textit{$r^p$-weighted},  quasi-conformal  vectorfields  defined in the non-causal,   far $r$ region of spacetime\footnote{These replace the scaling  and inverted time translation  vectorfields used in \cite{Kl-vectorfield}  or their  corresponding  deformations used in \cite{Ch-Kl}.  A recent improvement of the method  allowing one to derive higher order decay  can be found in \cite{AnArGa}.}.

   The first Energy-Morawetz  type  results  for scalar wave  equation \eqref{eq:WaveKerr}   in Schwarzschild, i.e. $a=0$,   are due to  Soffer-Blue  \cite{B-S1}, \cite{B-S2} and  Blue-Sterbenz \cite{B-St},    based on a modified version of the classical  Morawetz  integral energy decay  estimate.     Further developments appear in the works of Dafermos-Rodnianski, see  \cite{DaRo1},   \cite{Da-Ro3},  and   Marzuola-Metcalfe-Tataru-Tohaneanu   \cite{MaMeTaTo}. 
    The vectorfield method    can also be extended 	 to derive decay   for axially symmetric solutions   in Kerr, see    \cite{I-Kl1} and\footnote{In his Princeton PhD thesis Stogin   establishes a Morawetz estimate  even for the full subextremal   case $|a|<m$.} \cite{St}, but it is  known to fail   for  general  solutions  in Kerr, see  Alinhac \cite{Al}.

In the absence of axial symmetry   the derivation of    an    Energy-Morawetz estimate   in    $\KK(a,m)$ for  $|a/m|\ll 1 $ requires  a more refined  analysis involving     both the vectorfield method and either  micro-local methods or mode decompositions.   The first   full quantitative  decay\footnote{See also   \cite{DaRo2}  for the first proof of  boundedness of solutions, based on mode decompositions.}   result, based on micro-local analysis techniques, is due to Tataru-Tohaneanu  \cite{Ta-Toh}.  The   derivation of    such an estimate  in the full sub-extremal case $|a|<m$  is even more subtle  and   was  achieved  by Dafermos-Rodnianski-Shlapentokh-Rothman \cite{mDiRySR2014} by  combining the  vectorfield method with a  full  separation of variables  approach.    A  purely physical   space proof  of  the  Energy-Morawetz estimate   for small  $|a/m|$,  which avoids both  micro-local analysis  and mode  decompositions, was pioneered   by  Andersson-Blue in   \cite{A-B}.     Their method, which  extends the classical vectorfield method to  include    second order  operators (in this case  the Carter operator, see section
      \ref{section:Kerr}),      has the usual advantages of  the classical  vectorfield method, i.e it is robust with respect to perturbations.    It is for this reasons  that we  rely on it  in the proof of Theorem  \ref{MainThm-firstversion}, more precisely  in part II of  \cite{GKS-2022}.

       
\subsection{Linear  stability  of Schwarzschild}  
\lab{section:LinSchw}
  
  
 A  first  quantitative  proof   of the  linear stability of Schwarzschild  spacetime  was established\footnote{A  somewhat weaker version of  linear stability of Schwarzschild    was   subsequently   proved   in  \cite{HKW} by using the original, direct,   Regge-Wheeler, Zerilli  approach combined with the vectorfield method and    adapted gauge  choices. See also \cite{Johnson} for an alternate  approach  of  linear stability of Schwarzschild using wave coordinates.}  by  Dafermos-Holzegel-Rodnianski  (DHR) in \cite{D-H-R}.  
   Notable in their analysis  is  
  the treatment of  the Teukolsky equation  in  a fixed Schwarzschild background.   While  the Teukolsky equation    is separable, and  amenable to mode analysis, it  is not  variational and  thus cannot  be treated    directly  by      energy  type estimates. As mentioned earlier 
   in section \ref{section:Modestability},
    Chandrasekhar   was able to relate  the Teukolsky equation to the    Regge-Wheeler (RW)  equation, which is both  variational and coercive (the potential  $V$  has a  favorable sign).
    In   \cite{D-H-R}  the authors  rely   on a physical space  version of   the  Chandrasekhar transformation.  
    Once  decay  estimates  for   the RW  equation    have   been established, based on the technology developed earlier for the scalar wave equation in Schwarzschild,    the authors  recover   the  expected boundedness and  decay for  solutions to  the original   Teukolsky equation.  
    
      The remaining work in \cite{D-H-R}   is to       derive  similar control   for the other curvature  components and  the   linearized  Ricci coefficients associated to the double null foliation.  This   last step  requires carefully chosen gauge conditions, which the authors   make  within the framework of  a double null   foliation,  initialized  both  on the  initial hypersurface and the  background Schwarzschild horizon\footnote{The authors use a a scalar condition for the linearized lapse along the event horizon (part of what the authors call future normalized gauge), itself initialized from initial  data, see (212) and (214) in \cite{D-H-R}.}. This gauge fixing  from initial data leads to sub-optimal decay estimates  for some  of the metric coefficients\footnote{See (250)--(252)  and  (254) in \cite{D-H-R}.}   and is  thus inapplicable to the nonlinear case. This deficiency  was  fixed in  the PhD thesis of E. Giorgi,  in the context of the  linear stability of  Reissner-Nordstr\"om, see     \cite{Giorgi-T},  by relying on  a linearized version of the GCM construction  in \cite{KS}.

  
\subsection{Linear stability of  Kerr for small angular momentum}
\lab{section:LinstabKerr}


The first  breakthrough       result on the linear stability of Kerr,      for $|a|/m\ll 1$,  is due  to  Ma \cite{Ma}, see also  \cite{D-H-R-Kerr}. Both results  are  based on a generalization of the Chandrasekhar transformation to Kerr which   takes the Teukolsky equations, verified by the extreme  curvature components,    to   generalized versions of the  Regge-Wheeler (gRW)  equation.   
 Relying on separation of variables  and vectorfield techniques,  similar to those developed   for  the scalar wave equation in slowly rotating  Kerr, the authors  derive Energy-Morawetz and $r^p$ estimates    for the  solutions  of the gRW  equations. These  results   were  recently partially\footnote{The analysis  of  \cite{Y-R} is still  limited  to  modes. The authors   have  however announced a full proof    of  the result.}  extended to the full subextremal range in    \cite{Y-R}.
 
The first   stability results for the full linearized   Einstein vacuum equations   near   $\KK(a, m)$,   for $|a|/m\ll 1$,  appeared  in    \cite{ABBMa2019}  and  \cite{HHV}.  The first paper,  based on the   GHP formalism\footnote{An adapted  spinorial version  of the    NP formalism.}, see   \cite{GHP},    builds  on the results of  \cite{Ma}  while the second paper is based on an adapted  version of  the   metric formalism and builds on the seminal work of the authors on Kerr-de Sitter \cite{HVas}. Though the ultimate  relevance  of these  papers to  nonlinear stability  remains open,   they are  both  remarkable results  in so far as they deal with difficulties that looked insurmountable even ten years ago.


\subsection{Nonlinear model problems}  



\subsubsection{Nonlinear stability of Kerr-de Sitter}  


There is another important,  simplified,   nonlinear model problem which has drawn  attention
   in recent years, due mainly to the  striking achievement  of Hintz and Vasy \cite{HVas}. This is the problem of stability  of Kerr-de Sitter  concerning   
      the Einstein vacuum equation  with  a strictly positive cosmological 
constant
\bea
\lab{eq:Einstein-cosmological}
\R_{\a\b}  +\Lambda \g_{\a\b}=0, \qquad \La>0.
\eea

In their work, which   relies  in part   on the important  mode stability result     of   Kodama and  Ishibashi \cite{K-I},      Hintz and Vasy were able    prove the nonlinear   stability   of the stationary part  of  Kerr-de Sitter with small angular momentum, the first nonlinear stability  result of any nontrivial  stationary solutions   for the  Einstein equations\footnote{This is also   the first  general nonlinear  stability  result  in GR  establishing   asymptotic stability towards a family of solutions, i.e. full quantitative   convergence  to a final state  close,  but  different from the  initial one.}.
 It is important to note that,  despite  the   fact that,  formally,  the Einstein vacuum equation \eqref{eq:EVE-intro} is the limit\footnote{To pass  to the limit  requires one  to  understand  all global in time solutions of \eqref{eq:Einstein-cosmological} with $\La=1$, not only  those which are small perturbations of   Kerr-de Sitter,  treated by \cite{HVas}.} of \eqref{eq:Einstein-cosmological} as $\La\to 0$, the  global  behavior  of  the corresponding solutions   is   radically different\footnote{Major   differences between   formally close equations occur in many other   contexts. For example, the  incompressible Euler equations   are formally the limit of the Navier-Stokes equations as the viscosity tends to zero. Yet, at fixed  viscosity, the global properties of the Navier-Stokes equations are radically different   from that  of the  Euler equations.}.

  The main  simplification in the case of stationary solutions of  \eqref{eq:Einstein-cosmological} is that  the  expected decay   rates  of perturbations  near   Kerr-de Sitter     is exponential,   while in the case $\La=0$ the decay is  lower degree polynomial\footnote{While there is exponential decay in the stationary part treated in \cite{HVas}, note that lower degree polynomial decay is expected in connection to the  stability of the complementary causal region  (called   cosmological   or expanding)  of the full  Kerr-de Sitter  space, see e.g. \cite{Vo}.},  with various components of tensorial quantities  decaying at different  rates,  and the slowest decaying rate\footnote{Responsible  for carrying  gravitational waves at large distances so that they are  detectable.} being no better than $t^{-1}$.    The  Hintz-Vasy result was  recently revisited  in the work of A. Fang \cite{Fang2} \cite{Fang1} where he bridges the gap between the spectral methods of \cite{HVas} and the vectorfield methods.


\subsubsection{Nonlinear stability of Schwarzschild}  

 
The first  nonlinear stability result of the Schwarzschild space  was established  in  \cite{KS}.    In its simplest version, the  result  states the following.

\begin{theorem}[Klainerman-Szeftel \cite{KS}] 
\lab{MainThm-firstversion-KS}
The future globally hyperbolic   development  of  an   \emph{axially symmetric, polarized},  asymptotically  flat   initial data set, sufficiently close  (in a specified  topology)  to a Schwarzschild  initial data set   of  mass $m_0>0$,  has a complete    future null infinity  $\II^+$ and converges 
in  its causal past  $\JJ^{-1}(\II^{+})$  to another  nearby  Schwarzschild solution of mass $m_f$ close to $m_0$. 
 \end{theorem}

The restriction to  axial polarized perturbations  is the  simplest  assumption   which insures  that the final state is itself Schwarzschild  and thus avoids the additional complications  of  the Kerr stability problem.  We refer  the reader to  the introduction in \cite{KS} for a full discussion of the result.  

  The proof is based on  a construction based on GCM admissible spacetimes  similar to  that 
 briefly discussed  in section \ref{section:MainResult}  in the context of  slowly rotating Kerr.
 There are however several  important   simplifications to be noted:
\begin{itemize}
\item    The assumption of  polarization makes the  constructions of  the GCM  spheres $S_*$ and   spacelike hypersurface $\Si_*$ significantly simpler, see Chapter 9 in \cite{KS}, by comparison to  the general  case treated in  \cite{KS-GCM1},    \cite{KS-GCM2}  and  \cite{Shen}. 
\item The spacetime  has only two components   $\MM=\Mext\cup\Mint $ and the 
  null horizontal  structures, defined  on  each component,    are   integrable. 
 
 \item   As in the case of the scaler wave equation on Schwarzschild space the main    spin-2 Teukolsky   wave equations  can be treated by  a  vectorfield approach. This is no longer true 
    in   Kerr and even less so in  perturbations of Kerr.  
\end{itemize}
 
 \begin{figure}[h!]
\centering
\includegraphics[scale=0.35]{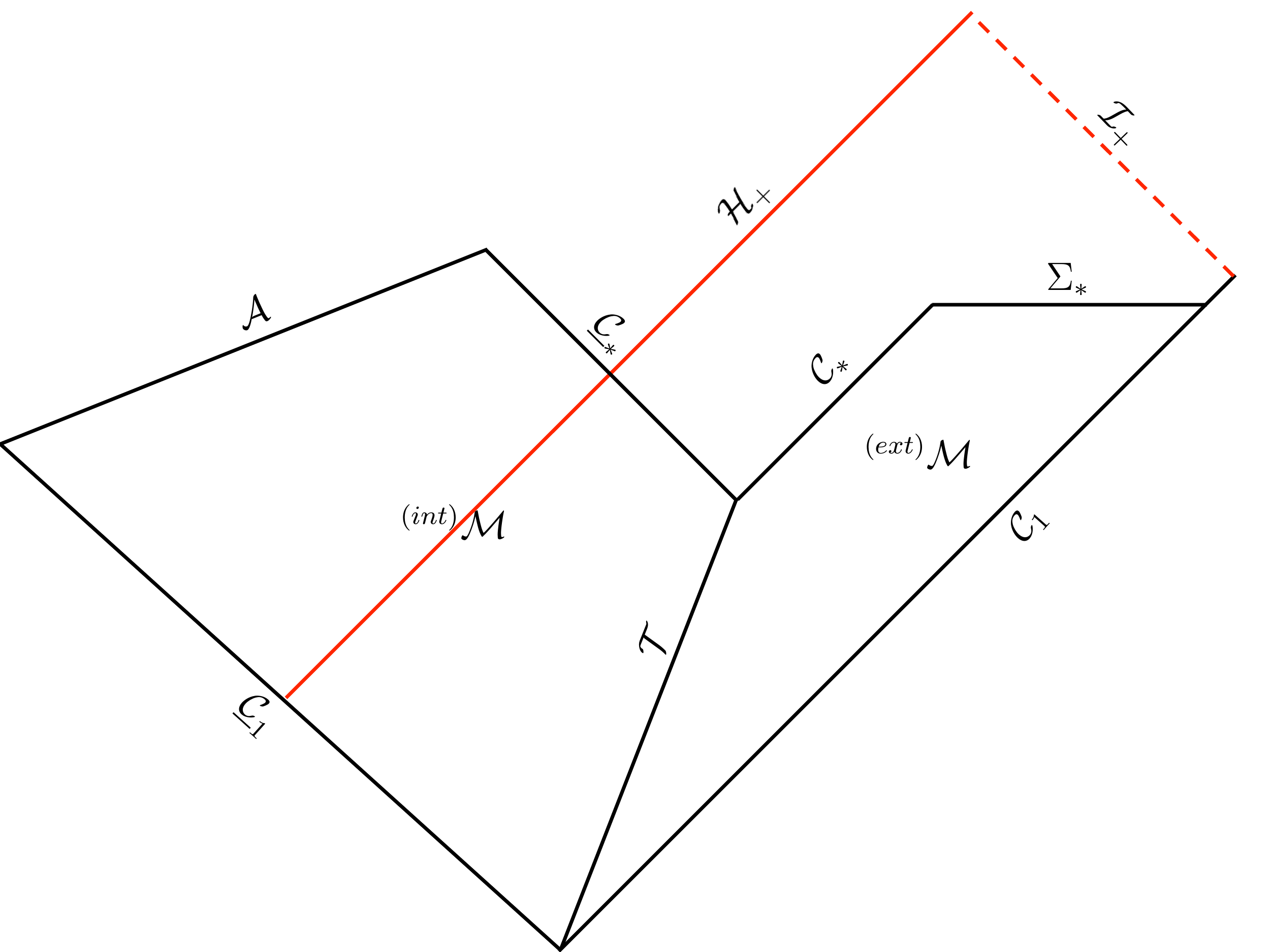}
\caption{\footnotesize{The GCM admissible space-time $\mathcal{M}$. By comparison to Figure \ref{fig1-introd1},  $\MM$  does not have $\Mtop$,  the past   boundaries 
$\CC_0\cup\CCb_0$ and future boundary   $\CCb\cup\CCb_*$  are null and   the horizontal structures (induced by  geodesic foliations)  are integrable.   As  in Theorem \ref{MainThm-firstversion}, the crucial GCM  sphere $S_*$ is  defined and constructed with no reference to the initial data. }}
\lab{fig:schwarzchildbootreg}
\end{figure}

Recently Dafermos-Holzegel-Rodnianski-Taylor   \cite{DHRT} have   extended\footnote{The novelty of \cite{DHRT},  compared to \cite{KS},  is the well preparation of  the initial data, based on an additional   three  dimensional modulation. Note however  that \cite{DHRT}  requires  substantially stronger   asymptotic conditions for the initial data compared to \cite{KS}.}  the result of \cite{KS}      by properly  preparing a co-dimension 3 subset of the initial data such that the final state is  still Schwarzschild. Like in   \cite{KS},  the starting point of  \cite{DHRT} is to anchor  the entire construction   on a  far away\footnote{That is  $r\gg u$, similar to the dominant in $r$ condition (3.3.4) of \cite{KS}.} GCM type  sphere $S_*$, in the sense of \cite{KS-GCM1}  \cite{KS-GCM2}, with no  direct reference to the initial data. It also  uses  the same definition of the  angular momentum as in (7.19) of \cite{KS-GCM2}. Finally, the spacetime in \cite{DHRT} is separated in an exterior region $\Mext$ and an interior region $\Mint$,   with the ingoing  foliation of $\Mint$,  initialized based on the information induced by $\Mext$,  as in \cite{KS}. We note, however, that \cite{DHRT} does not use the geodesic foliation of \cite{KS}, but instead both  $\Mint$ and $\Mext$ are foliated  by double null    foliations, and thus, the process of estimating the gauge dependent variables is somewhat different.

     
\section{Main ideas in the proof of Theorem \ref{MainThm-firstversion}}
\lab{section:Proof}



\subsection{The bootstrap  region}


 As mentioned in section \ref{section:Brief},  the proof   of Theorem  \ref{MainThm-firstversion}  is  centered   around    a  limiting  argument for    a family of  carefully constructed finite 
\textit{generally covariant modulated   (GCM) admissible  spacetimes} $\MM=\Mext\cup\Mtop\cup\Mint$.    As can be seen in Figure \ref{fig1-introd}  below, the  future boundary of the  spacetime is given by $\AA\cup\Sitop\cup \Si_*$ where  $\Si_*$ 
is a spacelike, generally covariant modulated (GCM) hypersurface, that is a hypersurface verifying  a  set of  crucial, well-specified,  geometric conditions,  essential to our proof of   convergence  to a final state. 
 
 \begin{figure}[h!]
 \centering
\includegraphics[scale=0.9]{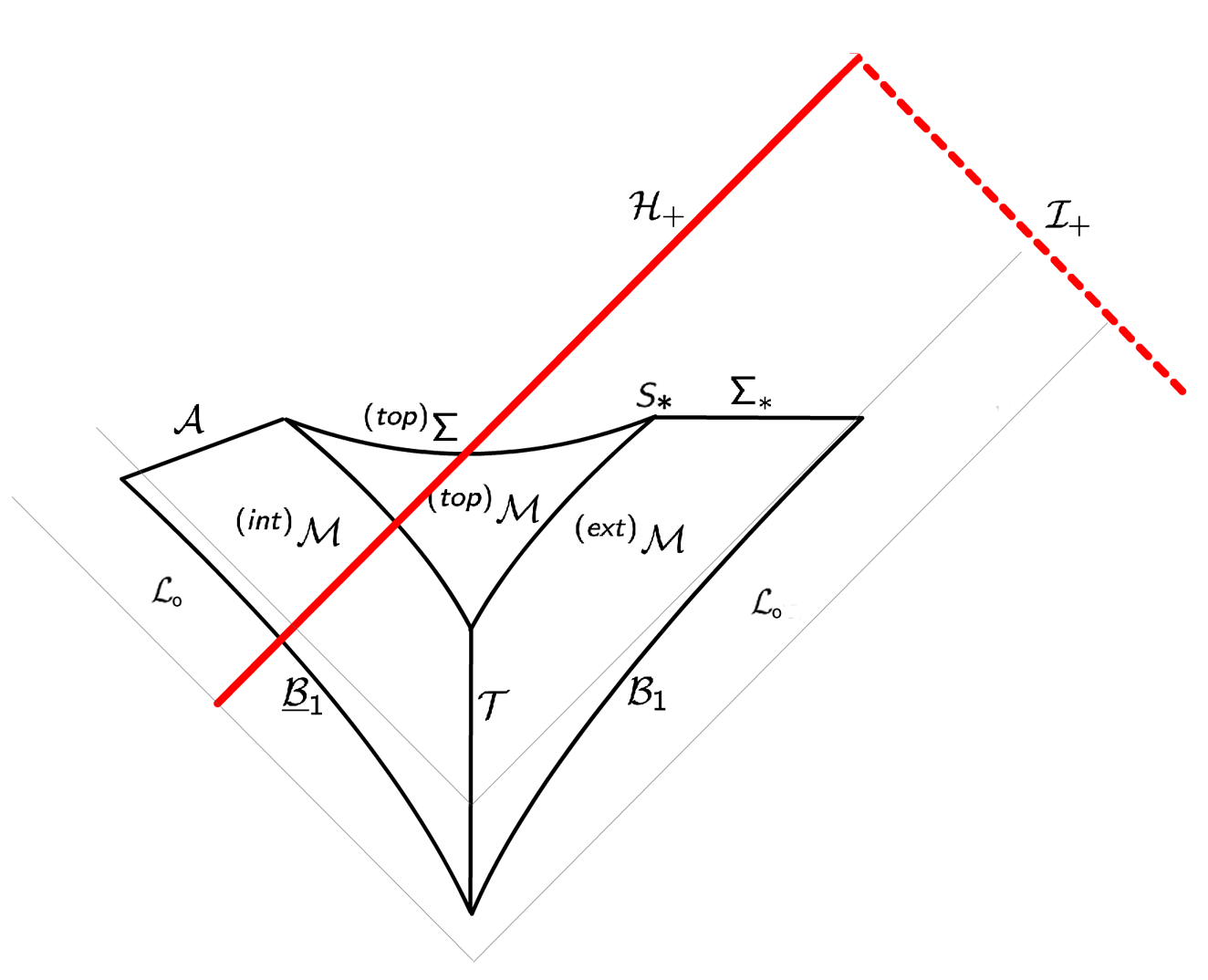}
\caption{\footnotesize{The Penrose diagram of  a finite GCM admissible    space-time $\MM=\Mext\cup\Mtop\cup\Mint$. The spacetime  is prescribed in the initial layer  $\LL_0$ and has $\AA\cup\Sitop\cup \Si_*$ as future boundary, with $\Si_*$  a spacelike ``generally covariant modulated (GCM)'' hypersurface. Its past boundary,  $\BB_1\cup\BBb_1$, is  itself part of the construction.  $\Mext$ is   initialized  by   the GCM hypersurface $\Si_*$ while  $\Mint$ is initialized  on  $\TT $ by the foliation induced by $\Mext$. The main inovation is the  GCM sphere $S_*$, defined  and constructed  with no reference  to the initial data prescribed   in   the  initial data layer $\LL_0$.  }}
\lab{fig1-introd} 
\end{figure} 
The capstone as well as  the most original part   of the entire construction is the  sphere $S_*$, the  future boundary  of $\Si_*$,  which verifies a  set of rigid, extrinsic and intrinsic,  conditions.  
 Once $\Si_*$ is specified the whole  GCM admissible spacetime $\MM$   is  determined by  a more  conventional  construction, based on  geometric transport type equations.
More precisely $\Mext$  can be determined  from $\Si_*$ by  a specified  outgoing  foliation terminating in the timelike boundary $\TT$,  $\Mint$ is determined from $\TT$ by a  specified  incoming 
 one, and $\Mtop$  is a complement of $\Mext\cup\Mint$   which makes $\MM$ a  
 causal domain\footnote{This is required because of the fact that, in our construction,  the  future  boundary of $\Mext\cup\Mint$ is not causal. By contrast,    in   \cite{KS}, $\MM=\Mext\cup \Mint$.}. The past boundary $\BB_1\cup\BBb_1$ of $\MM$, which is itself  to be constructed,  is included in the initial layer $\LL_0$ in which the  spacetime is assumed  to be known, i.e. a small perturbation of a Kerr solution.   The passage from  the initial data specified  on $\Si_0$ to  the initial layer  spacetime $\LL_0$ is  justified by D. Shen in \cite{Shen:Kerrext} by arguments similar to those of  \cite{KlNi1}-\cite{KlNi2},   based  on the mathematical methods  and techniques introduced in \cite{Ch-Kl}.

  Each of the spacetime regions $\Mext, \Mint,\Mtop $ come equipped  with  specific geometric  structure including specific   choices of  null   
    frames and functions such as $r,  u, \ub$. These are first defined on $\Si_*$ and then transported   to  $\Mext, \Mint,\Mtop $.


\subsection{Main intermediary results}


The proof of Theorem \ref{MainThm-firstversion} is divided in nine separate steps, Theorems M0--M8. These steps are  briefly described  below, see section 3.7 in \cite{KS:Kerr} for the precise statements:
\begin{enumerate}
 \item \textit{Theorem M0 (Control of the initial data in the bootstrap gauge)}.  The smallness of the initial perturbation  is   given  in the frame of the initial data layer $\LL_0$.  Theorem M0 transfers  this control to the bootstrap gauge in the initial data layer. 

 \item \textit{Theorems M1--M2 (Decay estimates for $\a$ (Theorem M1) and $\aa$ (Theorem M2))}.  This is achieved   using Teukolsky equations and a Chandrasekhar type transform in perturbations of Kerr. 

\item \textit{Theorems M3--M5 (Decay estimates for all curvature, connection and metric components)}. This is done making use of the GCM conditions on $\Si_*$ as well as    the control of $\a$ and $\aa$ established in Theorems M1 and M2. The proof proceeds in the following order:  
\begin{itemize}
\item Theorem M3 provides the  crucial  decay  estimates on   $\Si_*$, 

\item Theorem M4 provides  the  decay estimates on  ${}^{(ext)}\MM$,

\item Theorem M5 provides  the  decay estimates on  ${}^{(int)}\MM$ and ${}^{(top)}\MM$.
\end{itemize}

\item \textit{Theorems M6 (Existence of a bootstrap spacetime)}. This theorem shows that there exists a GCM admissible spacetime satisfying the bootstrap assumptions, hence initializing the bootstrap procedure.

\item \textit{Theorems M7 (Extension of the bootstrap region)}. This theorem shows the existence of a slightly larger GCM admissible spacetime satisfying estimates improving the bootstrap assumptions on decay.  

\item \textit{Theorem M8 (Control of the  top derivatives estimates)}.  This is based on an induction argument relative to the number of derivatives,  energy-Morawetz estimates and the Maxwell like character of the Bianchi identities. 
 \end{enumerate}

The paper  \cite{KS:Kerr} provides the proof of Theorem M0, Theorems M3 to M7, and half of Theorem M8 (on the control of Ricci coefficients and metric components). The proof of Theorems M1 and M2, and of the other half of Theorem M8 (on the control of curvature components), based  on   nonlinear wave equations techniques, are provided in  \cite{GKS-2022}. The construction of GCM spheres in   \cite{KS-GCM1} \cite{KS-GCM2}, and of GCM hypersurfaces in \cite{Shen} are used in the proof of Theorems M6 and M7 to construct respectively the terminal GCM sphere  $S_*$ and the last slice hypersurface $\Si_*$ from $S_*$.

   
 \subsection{Main new   ideas of the proof}
 \lab{section:Mainideas}
 
 
 Here is  a short description of the main new ideas in the proof  of  
 Theorem \ref{MainThm-firstversion} 
 and how they compare  with ideas  used in   other  nonlinear  results.

 
 \subsubsection{GCM admissible spacetimes}
 \lab{sec:GCMadmissiblespacetimes}
 
 
 \begin{itemize}
 \item As mentioned  already     the crucial concept  in the proof of Theorem  \ref{MainThm-firstversion}  is that of  a  GCM admissible spacetime,  whose  construction is  anchored  by the GCM sphere $S_*$ in Figure \ref{fig1-introd}.    GCM spheres\footnote{See the   discussion in the introductions to  \cite{KS-GCM1}, \cite{KS-GCM2}.}, are codimension 2 compact surfaces, unrelated   to the initial conditions,   on which  specific geometric quantities take Schwarzschildian values  (made possible by taking into account the full  general covariance of the Einstein vacuum equations). In addition to these  extrinsic conditions
 the sphere $S_*$ is endowed with a choice of ``effective\footnote{This  is meant to insure 
 the rigidity of the uniformization map, see \cite{KS-GCM2}.}  isothermal coordinates'',   $(\th, \vphi)$ verifying the following properties:

 \begin{itemize}
 \item  The metric on $S_*$  takes the form $ g= e^{2\phi} r^2\big( (d\th)^2+ \sin^2 \th (d\vphi)^2\big)$.
 
 \item   The integrals on $S_*$ of  the  $\ell=1$ modes\footnote{This is a natural generalization of $\ell=1$ spherical harmonics.}  
 $ J^{(0)} :=\cos\th$,   $ J^{(-)} :=\sin\th\sin\vphi$ and $ J^{(+)} :=\sin\th\cos\vphi $   vanish identically.
 \end{itemize}

\item Given the GCM sphere $S_*$ and the effective isothermal coordinates $(\th, \vphi)$ on it, our GCM procedure  allows us, in particular,  to define the   mass $m$, the angular momentum $a$ and a virtual axis of rotation which converge, in the limit,  to the final parameters $a_f, m_f$ and the  axis of rotation of the final Kerr\footnote{Previous  definitions of the angular momentum  in General Relativity  were given  in \cite{Rizzi}, \cite{Chen}, \cite{Chen2}, see also  \cite{Sz} for a comprehensive discussion of the subject.}.   We refer the reader to  section 7.2 in \cite{KS-GCM2} for our intrinsic definition of $a$ and of the virtual axis of symmetry on a GCM sphere.

 \item  The boundary $\Si_*$, called a GCM  hypersurface,  is initialized at $S_*$  and verifies additional conditions. 
  In the  polarized setting the  first such construction appears \cite{KS}.  The general case needed for our theorem is treated in   \cite{Shen}.   
 
\item   The concepts of         GCM   spheres   has appeared first in \cite{KS} in the context of polarized  symmetry.   The  construction of GCM spheres,  without any  symmetries, in realistic perturbations of Kerr, is treated in \cite{KS-GCM1}, \cite{KS-GCM2}\footnote{See also chapter 16 of \cite{DHRT} in the particular case of perturbations of Schwarzschild, where the same concept  appears instead under the name  ``teleological".}.   
 
 \item  The main novelty of  the GCM approach   is that it    relies on  gauge conditions initialized at a far away co-dimension $2$ sphere  $S_*$,  with no direct reference to the initial conditions. Previously known geometric constructions, such as  in \cite{Ch-Kl},   \cite{KlNi1} and \cite{KLR},     were based on codimension-1 foliations constructed on spacelike or null hypersurfaces and initialized    on the initial hypersurface\footnote{The first such  construction  appears in the proof of  the nonlinear stability of the Minkowski  space  \cite{Ch-Kl} where the ``inverse lapse foliation'' was  constructed on the ``last slice'',    initialized  at  spacelike  infinity $i^0$. Similar constructions, where the last slice is null rather than spacelike,  appear in \cite{KlNi1} and \cite{KLR}.}. Gauge conditions initialized from the future with no direct reference to the initial conditions, which was initiated in \cite{KS}, have since been used in other works, see \cite{Giorgi-T} \cite{graf} \cite{DHRT}.
  
 \item  The GCM  construction introduces the following new  important  conceptual   difficulty.  The  foliation on $\Si_*$, induced from the far away  sphere $S_*$, needs to be 
 connected, somehow,    to the initial  conditions (i.e. the initial layer $\LL_0$ in Figure \ref{fig1-introd}). This is achieved in  both \cite{KS} and  \cite{KS:Kerr} by transporting\footnote{That is, we transport the $\ell=1$ modes of some quantities from $S_*$ to $S_1$, see section 8.3.1 in \cite{KS:Kerr}.} the   sphere   $S_*$  to  a  sphere $S_1$ in the   the initial  layer        and  compare  it, using the rigidity properties 
of the GCM conditions, to a sphere of  the initial data layer. This  induces  a new foliation  of  the initial  layer  which  differs substantially from  the  original one, due to a shift  of the center of mass frame of the  final black hole, known  in the physics literature as a gravitational wave recoil\footnote{We refer the reader to  section 8.3  in \cite{KS:Kerr} for the details.}.
 \end{itemize}

 
 \subsubsection{Non integrability of the horizontal structure}
 \lab{sec:nonintegrabilityhorizontalstructure}
 
 
 As mentioned in section \ref{section:Kerr},  the canonical  horizontal structure   
  induced by the principal null directions $(e_3, e_4)$ in \eqref{eq:PrincipalNull}   of Kerr are non integrable.  The lack of integrability
  is    dealt with by    the  Newman-Penrose (NP)  formalism by  general null frames $(e_3, e_4, e_1, e_2)$, with  $e_1, e_2$ a specified  basis\footnote{Or rather  the complexified vectors $m=e_1+ ie_2$ and $\ov{m} = e_1-i e_2$.}   of the horizontal structure induced by the null pair   $(e_3, e_4)$.  It thus reduces all calculations to  equations  involving    the  Christoffel  symbols  of the frame, as scalar quantities. This  un-geometric feature of the formalism makes it difficult to use it in  the nonlinear setting of  the Kerr  stability  problem.  Indeed  complex calculations   depend    on higher  derivatives  of  all   connection coefficients of the NP  frame rather than only  those which are geometrically significant. 
   This  seriously affects  and complicates  the structure of non-linear corrections and makes  it difficult to avoid  artificial gauge type singularities\footnote{There are   no    smooth,  global   choices  of  a basis $(e_1, e_2)$. 
  The choice   \eqref{eq:canonicalHorizBasisKerr-intro} in Kerr, for example, is singular  at $\th=0, \pi$.}. This difficulty is avoided  in \cite{Ch-Kl} by working with  a  tensorial approach  adapted to   $S$-foliations, i.e. $\{e_3, e_4\}^\perp$ coincides, at every point, with the tangent space to $S$.  
  
  In our work we  extend,  with minimal  changes,   the tensorial approach  introduced in  \cite{Ch-Kl} to general  non-integrable  foliations.   
  The idea is very simple: we   define  Ricci coefficients $\chi, \chib,\eta, \etab, \ze, \xi, \xib,  \om, \omb$    exactly  as in \cite{Ch-Kl},  relative to an arbitrary   basis of  vectors $(e_a)_{a=1,2} $ of $\HH:=\{e_3, e_4\}^\perp$. In particular,  the null fundamental  forms  $\chi$ and $\chib$, are given by
    \beaa
\chib_{ab}&=\g(\D_ae_3, e_b),\qquad\quad\,\,\chi_{ab}=\g(\D_ae_4, e_b).
\eeaa
Due to the  lack of integrability of  $\HH$,  the null fundamental  forms  $\chi$ and $\chib$ are no longer symmetric.  They can be both decomposed as follows
\beaa
\chi_{ab}=\frac 1 2 \trch\de_{ab}+\frac1 2 \in_{ab}\atrch +\chih_{ab}, \qquad  \chib_{ab}=\frac 1 2 \trchb\de_{ab}+\frac1 2 \in_{ab} \atrchb +\chibh_{ab},
\eeaa
   where  the new scalars $\atrch$, $\atrchb$ measure the lack of integrability  of  the horizontal structure. The null curvature components are also defined 
 as in \cite{Ch-Kl}, 
\beaa
\a_{ab}=\R_{a4b4},\,\,\,\, \b_a=\frac 12 \R_{a434}, \,\,\,\,    \bb_a=\frac 1 2 \R_{a334},  \,\,\,\,  \aa_{ab}=\R_{a3b3}, \,\,\,\,  \rho=\frac 1 4 \R_{3434}, \,\,\,\,  \rhod=\frac 1 4\dual \R_{3434}.
\eeaa
   The  null structure and null Bianchi equations can then be derived  as  in the integrable case, see chapter 7 in  \cite{Ch-Kl}. The only new features are the presence of the scalars $\atrch, \atrchb $ in the  equations. Finally we note that the equations acquire additional simplicity if we pass to complex notations\footnote{The dual here is taken with respect to the antisymmetric  horizontal 2-tensor $\in_{ab}$.}, 
   \bea\lab{eq:defcomplexquantitieshorizontalstruct}
   \bsplit
A &:= \a+i\dual\a, \quad B:=\b+i\dual\b, \quad\, P:=\rho+i\dual\rho,\quad \Bb:=\bb+i\dual\bb, \quad \Ab:=\aa+i\dual\aa,
\\
 X&:= \chi+i\dual\chi, \quad \Xb:=\chib+i\dual\chib, \quad H:=\eta+i\dual \eta, \quad \Hb:=\etab+i\dual \etab, \quad 
 Z :=\ze+i\dual\ze.
\end{split}
\eea

 
\subsubsection{Frame transformations and choice of frames}


Given an arbitrary perturbation of Kerr,  there is no a-priori  reason to prefer an horizontal structure  to any other one obtained from the first by another perturbation of the same size.
  It is thus   essential  that we  consider  all possible frame transformations from  one   horizontal structure  $(e_4, e_3, \HH)$ to another one  $(e_4', e_3',  \HH')$  together with the transformation  formulas $\Ga\to \Ga'$, $R\to R'$  they generate. The most general  transformation formulas  between two null frames is given in Lemma  3.1  of \cite{KS-GCM1}. It depends on two  horizontal  $1$-forms  $f, \fb$ and a  real scalar function $\la$ and is given by
 \bea
 \lab{General-frametransformation-intro}
 \bsplit
  e_4'&=\la\left(e_4 + f^b  e_b +\frac 1 4 |f|^2  e_3\right),\\
  e_a'&= \left(\de_a^b +\frac{1}{2}\fb_af^b\right) e_b +\frac 1 2  \fb_a  e_4 +\left(\frac 1 2 f_a +\frac{1}{8}|f|^2\fb_a\right)   e_3,\\
 e_3'&=\la^{-1}\left( \left(1+\frac{1}{2}f\c\fb  +\frac{1}{16} |f|^2  |\fb|^2\right) e_3 + \left(\fb^b+\frac 1 4 |\fb|^2f^b\right) e_b  + \frac 1 4 |\fb|^2 e_4 \right).
 \end{split}
 \eea
The  very important   transformation formulas  $\Ga\to \Ga'$, $R\to R'$ are given in Proposition 3.3 of \cite{KS-GCM1}.

The case of $ \KK(a, m)$, $a\neq 0$   presents  an interesting  new  feature which can be described as follows:
\begin{itemize}
\item To capture  the  simplicity induced by the principle null directions in Kerr it is natural  to work with  non-integrable frames.
We do in fact  define all our main quantities  relative to  frames for which all quantities which vanish in Kerr are of  the size of the perturbation.

\item A   crucial  aspect of all    important  results in GR, based on   integrable  $S$- foliations,  is  that one can rely on elliptic Hodge theory  on each $2$-surface  $S$. This is  no longer possible in  context where  our main quantities  and the basic  equations  they verify are  defined  relative to non integrable frames.  In our work   we deal with this problem by passing back and forth, whenever needed,  
 from the main non-integrable frame  to a  well chosen adapted  integrable frame,   according to the transformation formulas mentioned above.
\end{itemize}

 
\subsubsection{Renormalization procedure and the canonical complex 1-form $\Jk$}
\lab{sec:renormalizationprocedure}


We first notice that  our main complex quantities introduced in \eqref{eq:defcomplexquantitieshorizontalstruct}
 take a particularly simple form in the principal null frame \eqref{eq:PrincipalNull} of Kerr:
\bea\lab{eq:particularsimpleforminprincipalnullKerr}
\bsplit
A &=\Ab=B=\Bb=0, \qquad P=-\frac{2m}{q^3},\\
\Xh &=\Xbh=0,\qquad \tr X=\frac{2}{q}\frac{\De}{|q|^2}, \qquad \tr\Xb=-\frac{2}{\ov{q}},\\
Z &=\frac{aq}{|q|^2}\Jk, \qquad H=\frac{aq}{|q|^2}\Jk, \qquad \Hb=-\frac{a\ov{q}}{|q|^2}\Jk,
\end{split}
\eea
where $q=r+ia\cos\th$, and where the regular\footnote{Note that $\Jk$ is regular including at $\th=0, \pi$.} complex 1-form $\Jk$ is given by 
 \bea
 \lab{eq:DefineJk}
 \Jk_1=\frac{i\sin\th }{|q|}, \qquad \Jk_2=\frac{\sin\th }{|q|},
 \eea 
see sections 2.4.2 and 2.4.3 in \cite{KS:Kerr}. In particular, the following holds for the complexified horizontal tensors of \eqref{eq:defcomplexquantitieshorizontalstruct} in the principal null frame \eqref{eq:PrincipalNull} of Kerr:
\begin{itemize}
\item the complex scalars $P$, $\tr X$ and $\tr\Xb$ are functions of $r$ and $\cos\th$,

\item the non vanishing complex 1-forms $H$, $\Hb$ and $Z$ consist of functions of $r$ and $\cos\th$ multiplied by $\Jk$,

\item the traceless symmetric complex 2-tensors $A$, $\Ab$, $\Xh$ and $\Xbh$ vanish identically.
\end{itemize}

Based on that observation, for a given horizontal structure perturbing the one of Kerr,  we can define a renormalization  procedure    by which, once we have\footnote{The constants $m$ and $a$ are computed on our GCM sphere $S_*$, see section \ref{sec:GCMadmissiblespacetimes}. $r$, $\th$ and $\Jk$ are chosen on $S_*$, transported to $\Si_*$ and then to $\MM$. The horizontal structure is also defined first on $\Si_*$ and then transported  to $\MM$.} suitable  constants $(a, m)$, suitable scalar functions $(r, \th)$,  and a  suitable   complex 1-form $\Jk$,  
and after  subtracting  the corresponding  values  in Kerr computed from $(a, m, r, \th, \Jk)$ for   all the  Ricci and  curvature  coefficients,  we  obtain   quantities  which are  first order in the perturbation.

  More precisely, once $(a,m)$,  $(r, \th)$ and $\Jk$   have been chosen, we renormalize the quantities in \eqref{eq:defcomplexquantitieshorizontalstruct} that do not vanish in Kerr as follows\footnote{The renormalization is written here in the case of a null pair $(e_3, e_4)$ with an ingoing normalization.}:
\bea\lab{eq:renormalizedquantities}
\bsplit
\widecheck{P}&:=P+\frac{2m}{q^3},\qquad \widecheck{\tr X}:=\tr X-\frac{2}{q}\frac{\De}{|q|^2}, \qquad \widecheck{\tr \Xb}:=\tr\Xb+\frac{2}{\ov{q}},\\
\widecheck{Z} &:=Z-\frac{aq}{|q|^2}\Jk, \qquad \widecheck{H}:=H-\frac{aq}{|q|^2}\Jk, \qquad \widecheck{\Hb}:=\Hb+\frac{a\ov{q}}{|q|^2}\Jk.
\end{split}
\eea


\subsubsection{Principal Geodesic and Principal Temporal structures}


In addition to the GCM gauge conditions on $\Si_*$, we need to construct a gauge on $\MM$ which relates the non integrable horizontal structure to the scalars $(r, \th)$ and the complex 1-form $\Jk$.   Two such gauges were introduced in \cite{KS:Kerr}:
\begin{itemize}
\item  Principal Geodesic (PG) structure, which is a generalization of the geodesic foliation to non-integrable horizontal structures,

\item  Principal Temporal (PT) structure, which favors transport equations along a null direction.
\end{itemize}

The PG structure is well suited for decay estimates, but  fails to be well posed.  Indeed,    due to the lack of  integrability  of the horizontal structure,  we cannot  control  the null structure equations\footnote{In integrable situation, like in the case of $S$-foliations,
the  Hodge systems on   the leaves of the $S$-foliation allows us to avoid the loss.} without a loss of derivative.
 The  PT structure, on the other hand,  is designed  so that the  loss  of derivatives    in the null structure equations, in the incoming or outgoing direction, is completely avoided. Note however that the PT structure  is not well suited to the derivation of decay estimates on $\Mext$ where $r$ can take arbitrary large values. In  \cite{KS:Kerr} we   work with both  gauge conditions, depending 
 on   the goal  we want to achieve, and rely on the transformation formulas  \eqref{General-frametransformation-intro} to pass from  one to the other.

In the outgoing normalization both  the outgoing  PG and PT structures consist  of  a choice   $(e_3, e_4, \HH)$, with $e_4$ null geodesic, together
 with a scalar functions $r, \th$ and a complex $1$-form $\Jk$  such that $e_4(r)=1$, $e_4(\th)=0 $, $\nab_4(q\Jk)=0$ where $q=r+ia\cos\th$.
 In addition: 
  \begin{enumerate}
  \item  In a  PG structure   the gradient  of $r$,  given by  $N=\g^{\a\b} \pr_\b r \pr_\a$,  is perpendicular  to  $\HH$,
  
  \item   In a    PT  structure  $\Hb=-\frac{a\ov{q}}{|q|^2}\Jk$, i.e. $\Hbc=0$ in view of \eqref{eq:renormalizedquantities}.
 \end{enumerate}
 
 A similar definition  of incoming  PG and PT  structures  is obtained  by interchanging the roles  of $e_3, e_4$.
  Note that both structures still  need to be initialized. The  outgoing  PG and PT structures of $\Mext$ are both initialized  on $\Si_*$ from the GCM frame of $\Si_*$,  while  the ingoing  PT structures of $\Mint$ and $\Mtop$ are initialized on the the timelike hypersurface $\TT$, see Figure \ref{fig1-introd}, using  the  data induced by the outgoing structures.


\subsubsection{Control of the extreme curvature components $A, \Ab$} 


It was already observed by Teukolsky that, in linear theory,   the extreme components of the curvature  are both  gauge invariant and   verify  decoupled  wave equations\footnote{See discussion in section  \ref{sect:LinearStab}.}.   In  our nonlinear context   this translates  to the statement that   the horizontal 2-tensors  $A, \Ab$, defined  relative to an $O(\ep)$ perturbation   of the principal  frame of  Kerr,  are $O(\ep^2)$-invariant, relative
  to  $O(\ep)$ frame transformations\footnote{This means that $f, \fb$, $\la-1$ are $O(\ep)$ in the transformation formulas \eqref{General-frametransformation-intro}.},  and  verify  tensorial wave equations of the form 
\bea
 \squared_2A  +L[A]&=& \err(\Gac, \Rc),\qquad   \squared_2\Ab   +\und{L}\, [\Ab ]= \und{\err}(\Gac, \Rc).
\eea
Here $\squared_2$ denotes the wave operator on horizontal symmetric traceless 2-tensors,   $L$ and $\und{L}$ are linear first order operators  and    $\Gac, \Rc$ denote the linearized Ricci  and curvature coefficients.  The error terms 
$ \err(\Gac, \Rc), \, \underline{ \err}(\Gac, \Rc)$ are nonlinear  expressions  in  $\Gac, \Rc$.      

     In linear theory, i.e.
when $\g$ is the Kerr metric and  the error terms  are not present,    these equations have been treated 
by  \cite{D-H-R} in Schwarzschild\footnote{See discussion in section \ref{section:LinSchw}.}  and by  \cite{Ma} and  \cite{D-H-R-Kerr} in slowly rotating\footnote{See discussion in section  \ref{section:LinstabKerr}.} Kerr, i.e. $|a|/m\ll 1$.   
   More precisely both results   derive realistic decay estimates for  $A, \Ab$. 
  The methods  are however  not    robust. Indeed,  a crucial ingredient in the proof,
    the   Energy-Morawetz estimates,  is based on  separation of variables. The control of $A$ and $\Ab$ in perturbations of Kerr in \cite{GKS-2022} contains the following new features:
   \begin{itemize}
   \item \textit{Derivation of the gRW equation}. The derivation of the      generalized Regge-Wheeler equations in Kerr, in  \cite{Ma} and  \cite{D-H-R-Kerr},  is done starting with  the  complex, scalar,  Teukolsky    equations, derived via  the  NP, or GHP   formalism,  
      by applying  a Chandrasekhar type transformation.  In part I of \cite{GKS-2022}  we  extend their derivation, using 
     our non-integrable horizontal formalism,  to perturbations of Kerr.  By contrast with    \cite{Ma}, \cite{D-H-R-Kerr}, we   derive  gRW equations for  the horizontal symmetric traceless 2-tensors\footnote{Derived from $A, \Ab$, see Definition 5.2.2 and 5.3.3 in \cite{GKS-2022}.}   $\qf, \qfb$,       rather than for complex scalars.   The main difficulty here is to make sure that the   non-linear  error terms  verify  a favorable structure.       
     
\item \textit{Nonlinear  error terms}.   The control of the nonlinear terms and their associated null structure was already understood in perturbations of Schwarzschild in \cite{KS} and is extended to perturbations of Kerr in \cite{GKS-2022}. 
\item \textit{Energy-Morawetz}. To derive energy-morawetz   estimates for $A, \Ab$ in  Part II of \cite{GKS-2022}
we  vastly extend the   pioneering   idea   of Andersson  and Blue \cite{A-B}, based on commutations with  $\T, \Z$ and the second order Carter operator $\CC$,   developed in the context of the scalar wave equation in slowly rotating Kerr,  to treat  our    tensorial  Teukolsky and gRW equations in  perturbations of Kerr.

   \end{itemize}


\subsubsection{Comment on the full  sub-extremal  range}


Though the full sub-extremal range  $|a|<m$  remains open  we remark that  a large part of our work 
does not require  the smallness of $|a|/m$.  This is the case 
 for  \cite{KS-GCM1} \cite{KS-GCM2} \cite{Shen} and \cite{KS:Kerr}. In fact  
 the   smallness  assumption is only needed in \cite{GKS-2022}, mostly  in the derivation of the main 
 energy-morawetz estimates in parts II and III.

\end{document}